\documentclass[3p]{elsarticle}

\usepackage{amsmath}
\usepackage{amsthm}
\usepackage{amssymb}
\usepackage{newtxtext,newtxmath}
\usepackage{graphicx}
\usepackage{hyperref}
\usepackage{booktabs}
\usepackage{enumitem}
\usepackage{subcaption}

\newtheorem{theorem}{Theorem}

\theoremstyle{definition}

\theoremstyle{remark}
\newtheorem{remark}{Remark}

\makeatletter
\let\ps@pprintTitle\ps@plain
\makeatother

\begin{document}

\begin{frontmatter}

\title{To Lift or Not to Lift: \\ Operator Inference for Nonlinear PDEs with Boundary Control}

\author[inst1]{Zhuoyuan Wang\fnref{merl}}
\author[inst1]{Yorie Nakahira}
\author[inst2]{Wael H. Ali}
\author[inst2]{Saviz Mowlavi\corref{cor1}}

\fntext[merl]{Work partially completed during an internship at Mitsubishi Electric Research Laboratories.}
\cortext[cor1]{Corresponding author. Email: \texttt{mowlavi@merl.com}}

\affiliation[inst1]{
    organization={Department of Electrical and Computer Engineering, Carnegie Mellon University},
    addressline={5000 Forbes Avenue},
    city={Pittsburgh},
    postcode={15213},
    state={PA},
    country={USA}
}

\affiliation[inst2]{
    organization={Mitsubishi Electric Research Laboratories},
    addressline={201 Broadway},
    city={Cambridge},
    postcode={02139},
    state={MA},
    country={USA}
}

\begin{abstract}

Operator inference (OpInf) provides a nonintrusive framework for learning projection-based reduced-order models directly from data, but its standard formulations do not explicitly account for the fact that boundary actuation changes the admissible solution space of a PDE. This paper develops and evaluates a boundary-consistent, lifting-based OpInf formulation for nonlinear PDEs with time-varying inhomogeneous boundary conditions. The proposed approach decomposes the state into a known lifting component that satisfies the prescribed boundary conditions and a homogeneous component represented in a POD basis. Substituting this decomposition into the governing dynamics reveals a boundary-consistent reduced model structure in which the learned dynamics includes bilinear state-input and quadratic input-input interactions for quadratically nonlinear PDEs. The resulting operators are inferred from data, preserving the nonintrusive character of OpInf while guaranteeing boundary condition satisfaction by construction. 
Numerical experiments on boundary-actuated heat, Burgers, and lid-driven cavity problems show that the lifted formulation improves boundary compliance and can enhance prediction accuracy, solver stability, and computational efficiency, although these gains vary across problems and reduced model dimensions. 
Overall, this study establishes lifting as a principled extension of OpInf for boundary-actuated PDEs, clarifies when and how this formulation improves performance, and indicates that standard OpInf can remain an effective alternative in certain regimes.

\end{abstract}


\end{frontmatter}

\section{Introduction}
\label{sec:introduction}

High-fidelity simulations of time-dependent partial differential equations (PDEs) are often computationally expensive, especially in repeated query tasks such as control, optimization, design, and uncertainty quantification. Projection-based Galerkin reduced-order models (ROMs) address this challenge by approximating the full-order state in a low-dimensional subspace, commonly obtained using proper orthogonal decomposition (POD), and evolving only the corresponding modal coefficients. Such methods have been widely used in fluid dynamics, flow control, and parametrized PDEs~\cite{rowley2017model,benner2015survey,hesthaven2016certified}. Classical POD-Galerkin ROMs preserve a close connection to the governing equations, but require access to the full-order operators and discretization details. This requirement can be restrictive for legacy solvers, commercial simulation codes, or experimental data. Data-driven approaches such as dynamic mode decomposition and sparse identification have therefore become popular alternatives for learning reduced dynamics directly from data~\cite{schmid2010dynamic,brunton2016discovering,rudy2017data}.

Operator inference (OpInf) provides a structured nonintrusive alternative by combining projection-based model reduction with data-driven regression. Instead of projecting known full-order operators, OpInf prescribes a reduced model form motivated by the Galerkin-projected equations and infers the corresponding reduced operators from trajectory data~\cite{peherstorfer2016data,kramer2024learning}. Under appropriate data richness and consistency assumptions, OpInf operators converge to the corresponding intrusive reduced operators, establishing a formal connection between data-driven and projection-based ROMs~\cite{peherstorfer2016data,benner2020operator}. OpInf and reduced operator inference have also been developed for broad classes of nonlinear PDEs, further clarifying the connection between learned reduced operators and the underlying projected PDE dynamics~\cite{qian2022reduced}. Recent work has extended OpInf in several directions, including regularization for stability and long-time prediction~\cite{mcquarrie2021data,sawant2023physics}, parametric PDEs~\cite{mcquarrie2023nonintrusive}, Bayesian uncertainty quantification~\cite{guo2022bayesian}, incompressible flows~\cite{benner2020operator2}, noisy or limited data~\cite{uy2023active,uy2023operator}, localized models~\cite{geelen2022localized}, and nonlinear trial manifolds~\cite{geelen2023operator}. Other structure-preserving extensions incorporate Hamiltonian, gradient, Lagrangian, and constrained forms of the reduced dynamics~\cite{sharma2022hamiltonian,geng2024gradient,sharma2024lagrangian,gruber2023canonical,kim2025physically,filanova2025reduced}. These developments demonstrate that the accuracy and stability of OpInf depend critically on choosing reduced coordinates and regression features that are consistent with the underlying physics.

Despite these advances, the treatment of time-dependent inhomogeneous boundary conditions in OpInf remains less developed. In many reduced models with inputs, the control is appended as an ad hoc additive forcing term in the reduced dynamics. This treatment is natural for distributed source terms, but it is generally not physically consistent for boundary actuation. Boundary inputs alter the admissible solution space itself: a linear combination of modes satisfying homogeneous boundary conditions cannot represent nonzero boundary values, while a POD basis constructed directly from raw snapshots must simultaneously represent both boundary offsets and interior dynamics, and will in general violate the prescribed boundary values. For quadratically nonlinear PDEs, boundary actuation can also induce bilinear state-input and quadratic input-input interactions. In intrusive ROMs, nonhomogeneous boundary conditions are commonly handled through lifting~\cite{gunzburger2007reduced}. The state is decomposed into a known component satisfying the prescribed boundary data and a homogeneous remainder that can be approximated using POD modes satisfying homogeneous boundary conditions. This idea has been used for time-dependent PDEs with parameters in the boundary data~\cite{gunzburger2007reduced}, Galerkin models and model predictive control for boundary-controlled Burgers and parabolic PDEs~\cite{hakenberg2013model,ou2009model,xu2011sequential}, boundary-condition-independent thermal ROMs~\cite{raghupathy2009boundary}, and incompressible flow ROMs with nonhomogeneous or time-dependent boundary conditions~\cite{ravindran2000reduced,barbagallo2009closed,grassle2019pod,star2019novel,hijazi2023pod}. Related boundary control and optimal control ROMs have also studied varying boundary controls and flow control actuation strategies~\cite{strazzullo2020pod,strazzullo2023pod,strazzullo2020reduced,kasnakouglu2008control}. However, these approaches are primarily developed in intrusive or semi-intrusive Galerkin settings, where the governing equations, weak forms, and boundary operators are available for projection.

This paper develops a boundary-consistent OpInf framework for nonlinear PDEs with time-dependent boundary actuation, with the goal of clarifying how the classical lifting treatment of nonhomogeneous boundary conditions can be translated into the nonintrusive OpInf setting and how it compares with standard OpInf formulations. The proposed lifting-based formulation follows the philosophy used in intrusive ROMs for nonhomogeneous boundary conditions~\cite{gunzburger2007reduced,barbagallo2009closed,hakenberg2013model,star2019novel}, while being related in spirit to existing transformation-based OpInf methods~\cite{qian2020lift,swischuk2020learning}, which use lifting for the distinct purpose of exposing a polynomial ROM structure from non-polynomial PDEs. We first decompose the state into a known lifting component satisfying the inhomogeneous boundary conditions, and a homogeneous component satisfying homogeneous boundary conditions and represented in a POD subspace. By substituting the lifted decomposition into the governing PDE, we identify the reduced operator structure induced by boundary actuation. In particular, for quadratically nonlinear PDEs, the learned dynamics include bilinear state-input and quadratic input-input interactions. The resulting operators are inferred from data, preserving the nonintrusive nature of OpInf while retaining the boundary consistency of lifting-based Galerkin ROMs. This construction guarantees boundary condition satisfaction by design and yields a reduced representation in which the POD basis approximates the homogeneous interior dynamics. Our numerical results show that lifting enforces boundary compliance and can improve accuracy and stability in boundary-actuated problems. At the same time, we observe that standard OpInf formulations can still perform well without an explicit lifting design, especially when the primary objective is accurate trajectory prediction rather than exact boundary consistency. Thus, rather than suggesting that lifting is always necessary, this study clarifies what can be gained by lifting, what standard OpInf can already achieve, and when each formulation is justified.

The main contributions of this work are summarized as follows:
\begin{itemize}
\item We identify and address a gap between data-driven OpInf and classical reduced-order modeling for nonlinear PDEs with boundary control: while lifting is standard in intrusive ROMs for enforcing nonhomogeneous boundary conditions, its implications for the reduced coordinates, regression features, and learned operator structure in nonintrusive OpInf have not been systematically formulated.
\item We derive a boundary-consistent OpInf formulation for nonlinear PDEs with time-dependent inhomogeneous boundary conditions. The derivation shows how boundary lifting may induce bilinear state-input and quadratic input-input interactions, thereby providing a principled OpInf model structure for boundary-actuated systems with guaranteed satisfaction of the prescribed boundary conditions.
\item We demonstrate the method on boundary-actuated heat, Burgers, and lid-driven cavity problems, comparing the lifted formulation with standard OpInf alternatives and showing what can be gained in projection and prediction accuracy, boundary compliance, and stability. The results also clarify when standard OpInf remains effective without an explicit lifting design, thereby helping justify the choice between lifted and standard formulations.
\end{itemize}

The remainder of the paper is organized as follows. Section~\ref{sec:problem} introduces the problem formulation and reviews Galerkin projection and standard OpInf for PDEs with homogeneous boundary conditions. Section~\ref{sec:method} presents the proposed lifting-based boundary-consistent OpInf formulation and proves boundary condition satisfaction. Section~\ref{sec:experiments} evaluates the method on boundary-actuated heat and Burgers equations, as well as incompressible flow in a lid-driven cavity. The results are discussed in Section~\ref{sec:discussion}, and Section~\ref{sec:conclusion} concludes the paper.

\section{Problem Formulation and Background on Operator Inference}
\label{sec:problem}

Let $\Omega \subset \mathbb{R}^d$ be a bounded domain with Lipschitz continuous boundary $\partial\Omega$ and let $[0, T] \subset \mathbb{R}$ be a time interval of interest. Let $\mathcal{H}$ denote the Hilbert space of real-valued functions on $\Omega$ equipped with the standard $L^2$ inner product $\langle \cdot, \cdot \rangle$. We consider the general class of quadratically nonlinear partial differential equations (PDEs) subject to boundary control,
\begin{subequations}\label{eq:generic_pde}
\begin{align}
    \frac{\partial y}{\partial t} &= \mathcal{L}y + \mathcal{Q}(y, y) \quad \text{in } \Omega \times [0, T], \\
    \mathcal{B}y &= \mathcal{C} u(t) \quad \text{on } \partial \Omega \times [0, T], \label{eq:inhomogeneous_bc}
\end{align}
\end{subequations}
where $y(\cdot, t) \in \mathcal{H}$ is the system state, $\mathcal{L}$ is a linear operator, $\mathcal{Q}$ is a bilinear operator, $\mathcal{B}$ is a linear boundary operator, $u(t) \in \mathbb{R}^m$ are control inputs, assumed to be smooth, and $\mathcal C$ maps control inputs to boundary data \footnote{For example, $\mathcal C$ may prescribe boundary data through a finite set of actuator or boundary shape functions, $\mathcal C u(t) = \sum_{j=1}^m u_j(t) h_j(x)$, where $u_j(t)$ are scalar control amplitudes and $h_j$ are prescribed functions on the controlled portion of $\partial\Omega$.}.
We further assume throughout that the operators and function spaces are chosen so that the above initial boundary value problem is well posed. Contrary to most OpInf papers, we derive the structure of the reduced-order model in the spatially continuous setting in order to expose the role played by the time-varying boundary conditions.

Our goal is to construct a data-driven reduced-order model (ROM) of the PDE \eqref{eq:generic_pde} from high-fidelity trajectory data. Starting from the case of homogeneous BCs for which $u = 0$, we first review standard model-based Galerkin projection in Section~\ref{sec:galerkin} and data-driven operator inference in Section~\ref{sec:opinf}. Then, we list the challenges faced by these methods when handling PDEs with inhomogeneous BCs in Section~\ref{sec:challenges}.

\subsection{Model-Based Galerkin Projection}
\label{sec:galerkin}
Galerkin projection~\cite{rowley2017model} is a systematic method to obtain a ROM for a PDE with homogeneous BCs. The core idea is to approximate the PDE dynamics in a finite-dimensional subspace spanned by basis functions $\phi_1, \ldots, \phi_r \in \mathcal{H}$, chosen to be orthonormal with respect to $\langle\cdot, \cdot\rangle$ and to satisfy the homogeneous BCs of the PDE. While such basis functions can be selected in different ways, the proper orthogonal decomposition (POD) \cite{holmes2012turbulence} is a popular data-driven approach using state snapshots collected along one or multiple solution trajectories of the PDE. Once these basis functions are obtained, the state is approximated within the span of the POD basis as
\begin{equation}
\label{eq:pod_expansion}
    y_r(x, t)=\sum_{i=1}^r a_i(t) \phi_i(x),
\end{equation}
where $a_i(t) \in \mathbb{R}$ are time-dependent modal coefficients. The expansion also satisfies the homogeneous BCs due to the linearity of $\mathcal{B}$. The dynamics for the modal coefficients are then obtained by substituting this approximation into the PDE and enforcing that the residual error be orthogonal to the subspace spanned by the basis functions. That is,
\begin{equation}
\label{eq:residual_condition}
    \left\langle\frac{\partial y_r}{\partial t}-\mathcal{L} y_r - \mathcal{Q}(y_r, y_r), \phi_j\right\rangle=0, \quad j=1, \ldots, r.
\end{equation}
Substituting \eqref{eq:pod_expansion} into \eqref{eq:residual_condition} and leveraging the orthonormality of the basis functions yields a reduced system for the modal coefficients
\begin{equation}
\label{eq:galerkin_rom}
    \frac{d a}{d t}=A a+H(a \otimes a),
\end{equation}
where $a=\left[a_1, \ldots, a_r\right]^T \in \mathbb{R}^r$ is the reduced state, $\otimes$ denotes the Kronecker product, $A \in \mathbb{R}^{r \times r}$ is the reduced linear operator with entries
\begin{equation}
\label{eq:galerkin_A}
    A_{i j}=\left\langle\mathcal{L} \phi_j, \phi_i\right\rangle,
\end{equation}
and $H \in \mathbb{R}^{r \times r^2}$ is the reduced quadratic operator with entries
\begin{equation}
\label{eq:galerkin_H}
    H_{i m}=\left\langle \mathcal{Q}\left(\phi_j, \phi_k\right), \phi_i\right\rangle, \quad m=(j-1) r+k.
\end{equation}

During inference, the initial condition for the modal coefficients at time $t_0$ is obtained from $a_i(t_0) = \langle y(\cdot, t_0), \phi_i \rangle$ for $i = 1, ..., r$, where $y(\cdot, t_0)$ is the prescribed initial full-order state. From the reduced trajectory calculated by \eqref{eq:galerkin_rom}, the full-order state trajectory is reconstructed using the expansion \eqref{eq:pod_expansion}. 

Galerkin projection thereby transforms a PDE with homogeneous BCs into a ROM that preserves the intrinsic structure of the PDE. However, its derivation requires access to the PDE operators $\mathcal{L}$ and $\mathcal{Q}$ which may unknown or inaccessible.

\subsection{Data-Driven Operator Inference}
\label{sec:opinf}
Contrary to Galerkin projection, operator inference (OpInf)~\cite{peherstorfer2016data,kramer2024learning} seeks to construct a ROM in a fully data-driven manner.
Assume that we have collected a trajectory dataset, either through experiments or simulations of the PDE~\eqref{eq:generic_pde}, consisting of a time series of full-order states
\begin{equation}
\label{eq:training_dataset_opinf}
    \left\{y(\cdot, t_1), \ldots, y(\cdot, t_K)\right\} \subset \mathcal{H}.
\end{equation}
Although~\eqref{eq:training_dataset_opinf} is written for a single trajectory for notational convenience, in practice the dataset typically contains states from multiple trajectories corresponding to different initial conditions.

OpInf uses these data to infer a ROM that is physically consistent in the sense that it retains the reduced operator structure of the corresponding Galerkin-projected ROM.
Accordingly, it adopts the same state approximation~\eqref{eq:pod_expansion} as Galerkin projection, with the POD basis computed from the snapshots in \eqref{eq:training_dataset_opinf}, and then postulates a reduced system with a similar structure to the Galerkin-projected ROM for the modal coefficients,
\begin{equation}
    \frac{d a}{d t}=A a + H(a \otimes a) + c,
\label{eq:postulated_rom_homogeneous}
\end{equation}
where $a=\left[a_1, \ldots, a_r\right]^T \in \mathbb{R}^r$ is the reduced state, $A \in \mathbb{R}^{r \times r}$ is a reduced linear operator, $H \in \mathbb{R}^{r \times r^2}$ is a reduced quadratic operator, and $c \in \mathbb{R}^{r \times 1}$ is a constant offset vector often included in OpInf formulations to improve prediction accuracy. 

Rather than computing the operators $A$ and $H$ from the Galerkin projection formulas in~\eqref{eq:galerkin_A} and~\eqref{eq:galerkin_H}, OpInf infers them directly from the data in \eqref{eq:training_dataset_opinf}. To this end, the full-order snapshots are first projected onto the POD basis to obtain a time series of reduced states and their time derivatives
\begin{equation}
    \{a(t_1), \dots, a(t_K)\} \subset \mathbb{R}^r, \quad \{\dot{a}(t_1), \dots, \dot{a}(t_K)\} \subset \mathbb{R}^r,
\end{equation}
where the reduced states are computed as
\begin{equation}
    a(t_k)=\left[\left\langle y(\cdot, t_k), \phi_1\right\rangle, \dots, \left\langle y(\cdot, t_k), \phi_r\right\rangle\right]^T, \quad k=1, \ldots, K \label{eq:projected_states}
\end{equation}
and the corresponding derivatives $\dot{a}(t_k)$, $k=1, \ldots, K$, are typically approximated by applying finite differences to the time series of reduced states~\cite{kramer2024learning}. When temporal derivatives of the full-order state are available from the PDE solver, the reduced time derivatives may also be obtained by directly projecting the full-order derivatives onto the POD basis in the same manner as in~\eqref{eq:projected_states}. Finally, OpInf calculates the operators $A, H$, and $c$ such that the postulated ROM in \eqref{eq:postulated_rom_homogeneous} best fits the data in the least-squares sense, which can be formulated~\cite{qian2022reduced} as the linear regression problem
\begin{equation}
\label{eq:OpInfMinimization}
    \min _{A, H, c} \sum_{k=1}^K\left\|A a(t_k)+H(a(t_k) \otimes a(t_k))+c-\dot{a}(t_k)\right\|_F^2,
\end{equation}
where $\|\cdot\|_F$ denotes the Frobenius norm. Penalty terms involving the Frobenius norms of the inferred operators can be appended to the least-squares objective to improve well-posedness of the minimization problem and robustness of the learned ROM when facing limited or noisy training data~\cite{mcquarrie2021data,qian2022reduced,sawant2023physics,kramer2024learning}.

Thus, OpInf yields a data-driven ROM consistent with the mathematical structure of the PDE with homogeneous BCs, without requiring knowledge of the PDE operators $\mathcal{L}$ and $\mathcal{Q}$.

\subsection{Challenges in the Presence of Boundary Actuation}
\label{sec:challenges}
Let us now consider PDEs with inhomogeneous BCs resulting from nonzero control inputs $u$ in~\eqref{eq:generic_pde}. We seek a data-driven ROM that describes the dynamics of such PDEs in response to arbitrary sequences of control inputs $u(t)$, where the trajectory dataset now also contains a time series of control inputs $\left\{u(t_1), \ldots, u(t_K)\right\} \subset \mathbb{R}^m$. Since Galerkin projection does not readily apply to PDEs with inhomogeneous BCs, OpInf addresses this by retaining the same expansion \eqref{eq:pod_expansion} while appending an ad hoc control term to the ROM derived for the homogeneous case, yielding
\begin{equation}
\label{eq:opinf_adhoc}
    \frac{d a}{d t}=A a+H(a \otimes a)+B u+c,
\end{equation}
where $B \in \mathbb{R}^{r \times m}$ is a control input matrix. The least squares problem becomes
\begin{equation}
    \min _{A, H, B, c} \sum_{k=1}^K\left\|A a(t_k)+H(a(t_k) \otimes a(t_k))+B u(t_k)+c-\dot{a}(t_k)\right\|_F^2 .
\end{equation}
Because this approach relies on an ad hoc modification of the ROM, it creates two problems:
\begin{enumerate}
    \item The expansion \eqref{eq:pod_expansion} is unable to represent arbitrary boundary values set by the boundary condition $\mathcal{B} y=\mathcal{C} u(t)$ on $\partial \Omega$ for arbitrary $u(t)$. This will inevitably lead to larger errors near the boundaries in the ROM predictions.

    \item The modified ROM \eqref{eq:opinf_adhoc} is no longer guaranteed to preserve the mathematical structure of the governing PDE, since its form is not derived from a Galerkin projection. For example, in quadratically nonlinear PDEs, the control input may affect the dynamics through bilinear state-input interactions, yet such coupling terms are missing from the modified ROM.
\end{enumerate}
As a result, standard OpInf ROMs may lose accuracy when applied to PDEs with inhomogeneous boundary conditions. This motivates a data-driven formulation that preserves the operator structure induced by boundary actuation while ensuring consistency with the prescribed boundary conditions.

\section{Proposed Method: Lifted Operator Inference}
\label{sec:method}

In this section, we outline our method to construct a data-driven yet physically consistent ROM of a governing PDE with inhomogeneous and time-dependent BCs. 
Our method adapts the classical lifting treatment of inhomogeneous BCs in Galerkin ROMs to the OpInf setting, enforcing boundary compliance by design and giving rise to previously ignored terms in standard OpInf formulations of boundary-actuated quadratic PDEs.

Consider the PDE with boundary actuation described in~\eqref{eq:generic_pde}, where $u(t)=[u_1, \ldots, u_m]^T\in\mathbb{R}^m$. 
Assume that we have collected a trajectory dataset, either through simulations or experiments, of full-order states and control inputs
\begin{equation}
\label{eq:training_dataset_proposed}
    \left\{y(\cdot, t_1), \ldots, y(\cdot, t_K)\right\} \subset \mathcal{H}, \quad \left\{u(t_1), \ldots, u(t_K)\right\} \subset \mathbb{R}^m.
\end{equation}
Following the classical literature on Galerkin ROMs for inhomogeneous PDEs~\cite{hakenberg2013model,ravindran2000reduced,barbagallo2009closed}, our approach begins with the exact state decomposition
\begin{equation}
\label{eq:decomposition}
y(x,t)=\sum_{j=1}^m u_j(t)\bar{g}_j(x)+\bar{y}_H(x,t),
\end{equation}
where $\bar{g}_j\in\mathcal{H}$ are prescribed lifting functions and $\bar{y}_H(\cdot,t)\in\mathcal{H}$ is a remainder. The lifting functions are chosen to satisfy 
\begin{equation}
\mathcal{B}\bar{g}_j=\mathcal{C}e_j \quad \text{on }\partial\Omega,\quad j=1,\ldots,r.
\end{equation}
with $e_j\in\mathbb{R}^m$ denoting the $j$-th standard basis vector. In other words, each lifting function $\bar{g}_j$ satisfies the inhomogeneous BCs \eqref{eq:inhomogeneous_bc} for the particular input $u(t)=e_j$ whose $j$-th component equals one and all others are zero\footnote{Such lifting functions can be readily constructed for various types of BCs. For example, consider the one-dimensional domain $[0,1]$ with two-dimensional control input $u(t)=[u_1(t),u_2(t)]^T$. For Dirichlet BCs $y(0,t)=u_1(t)$ and $y(1,t)=u_2(t)$, one may choose $\bar{g}_1(x)=1-x$ and $\bar{g}_2(x)=x$. For the mixed BCs $y(0,t)=u_1(t)$ and $\partial y(1,t)/\partial x=u_2(t)$, one may choose $\bar{g}_1(x)=1+2x-x^2$ and $\bar{g}_2(x)=x$.}. This construction ensures that the remainder $\bar{y}_H$ satisfies homogeneous BCs at all times, that is $\mathcal{B}\bar{y}_H=0$, owing to the linearity of $\mathcal{B}$.

We then construct orthonormal basis functions $\phi_1, \ldots, \phi_r \in \mathcal{H}$ satisfying homogeneous BCs, for example through a POD of the homogeneous remainder snapshots
\begin{equation}
\label{eq:training_dataset_homogenized}
    \left\{\bar{y}_H(\cdot, t_1), \ldots, \bar{y}_H(\cdot, t_K)\right\} \subset \mathcal{H},
\end{equation}
where each remainder snapshot is computed from \eqref{eq:decomposition}. Using this basis, we approximate the full state as
\begin{equation}
\label{eq:reconstruction}
    y_r(x, t)=\sum_{j=1}^m u_j(t)g_j(x) + \sum_{i=1}^r a_i(t) \phi_i(x),
\end{equation}
{where $a_i(t) \in \mathbb{R}$ are time-dependent modal coefficients and $g_j\in\mathcal{H}$ are `projected' lifting functions defined as
\begin{equation}
\label{eq:projected_lifting}
    g_j=\bar{g}_j - \sum_{k=1}^r\langle \bar{g}_j,\phi_k\rangle\phi_k,\quad j=1,\ldots,m.
\end{equation}
These projected lifting functions satisfy the orthogonality property 
$\langle g_j,\phi_i\rangle=0$ for $i=1, \dots, r$ and $j=1, \dots, m$, which will simplify the structure of the modal dynamics derived in the next step. Moreover,} the following theorem shows that the state approximation \eqref{eq:reconstruction} is guaranteed to satisfy the inhomogeneous BCs.

\begin{theorem}
\label{thm:bc_satisfaction}
The lifted state approximation~\eqref{eq:reconstruction}, using POD basis functions $\phi_1,\ldots,\phi_r$ obtained from the remainder snapshots~\eqref{eq:training_dataset_homogenized}, is guaranteed to satisfy the inhomogeneous BCs~\eqref{eq:inhomogeneous_bc}. That is,
\begin{equation}
\label{eq:bc_satisfaction}
\mathcal{B}y_r=\mathcal{C}u(t) \quad \text{on }\partial\Omega.
\end{equation}
\end{theorem}

\begin{proof}
Since the POD basis is constructed from the homogeneous remainder snapshots in~\eqref{eq:training_dataset_homogenized}, each mode $\phi_i$ lies in their linear span~\cite{volkwein2013proper,ullmann2014pod}. By construction, these snapshots satisfy the homogeneous BCs; therefore, by the linearity of the boundary operator $\mathcal{B}$, so do the POD modes, i.e.,
\begin{equation}
\label{eq:modes_bcs}
\mathcal{B}\phi_i=0 \quad \text{on }\partial\Omega,\quad i=1,\ldots,r.
\end{equation}
{This ensures that the projected lifting functions $g_j$ retain the same BCs as their raw counterparts $\bar{g}_j$. Indeed,
\begin{equation}
\label{eq:projected_lifting_functions_bcs}
\mathcal{B}g_j=\mathcal{B}\bar{g}_j - \sum_{k=1}^r\langle \bar{g}_j,\phi_k\rangle \mathcal{B}\phi_k=\mathcal{C}e_j \quad \text{on }\partial\Omega,\quad j=1,\ldots,r,
\end{equation}}
where we have again leveraged the linearity of $\mathcal{B}$. Finally, using the latter with \eqref{eq:modes_bcs} and \eqref{eq:projected_lifting_functions_bcs} yields
\begin{equation}
\mathcal{B}y_r=\mathcal{B}\left[\sum_{j=1}^m u_j(t)g_j+\sum_{i=1}^r a_i(t)\phi_i\right]=\sum_{j=1}^m u_j(t)\mathcal{B}g_j=\sum_{j=1}^m u_j(t)\mathcal{C}e_j=\mathcal{C}u(t),\quad \text{on }\partial\Omega.
\end{equation}
\end{proof}

As in Galerkin projection, the structure of the reduced dynamics for the modal coefficients is obtained by substituting the state approximation \eqref{eq:reconstruction} into the PDE and enforcing that the resulting residual be orthogonal to the subspace spanned by the basis functions. That is, we require
\begin{equation}
\label{eq:residual_condition_inhomogeneous}
    \left\langle\frac{\partial y_r}{\partial t}-\mathcal{L}y_r-\mathcal{Q}(y_r, y_r), \phi_j\right\rangle=0, \quad j=1, \ldots, r.
\end{equation}
Substituting \eqref{eq:reconstruction} in \eqref{eq:residual_condition_inhomogeneous} and leveraging the orthonormality of the basis functions, as well as the orthogonality between the basis and the projected lifting functions, yields a reduced system for the modal coefficients taking the form
\begin{equation}
\label{eq:opinf_proposed}
\begin{aligned}
    \frac{d a}{d t}
    &= A \xi + H(\xi \otimes \xi) + c \\
    &= A_a a + H_{aa}(a\otimes a) + A_u u + c 
    + \underbrace{
    H_{au}(a\otimes u)
    + H_{ua}(u\otimes a)
    + H_{uu}(u\otimes u)
    }_{\text{additional terms in the lifted formulation}} .
\end{aligned}
\end{equation}
where $a=[a_1, \ldots, a_r]^T \in \mathbb{R}^r$ is the reduced state, $\xi=[a^T, u^T]^T \in \mathbb{R}^{r+m}$ is an augmented reduced state, $A \in \mathbb{R}^{r \times(r+m)}$ is a reduced linear operator, and $H \in \mathbb{R}^{r \times(r+m)^2}$ is a reduced quadratic operator. In the expanded form, $A_a$ and $A_u$ denote the linear state and control operators, respectively, while $H_{aa}$, $H_{au}$, $H_{ua}$, and $H_{uu}$ denote the quadratic state, state-control interaction, control-state interaction, and quadratic control operators, respectively. In addition, $c \in \mathbb{R}^{r \times 1}$ is a constant offset.

\begin{remark}
The augmented reduced state $\xi$ provides a compact notation for the modal dynamics \eqref{eq:opinf_proposed} and subsequent minimization problem \eqref{eq:opinf_minimization_proposed}. However, only the modal coefficients in $a$ are evolved by the ROM; the control input $u$ is prescribed externally.
\end{remark}
\begin{remark}
The expanded modal dynamics in \eqref{eq:opinf_proposed} reveals several input-dependent contributions that are missing from the ad hoc modal dynamics \eqref{eq:opinf_adhoc}: terms bilinear in $a$ and $u$ as well as terms quadratic in $u$, which are inherited from the quadratic PDE operator $\mathcal{Q}$.
\end{remark}
\begin{remark}
The absence of a term proportional to $\dot{u}$ in \eqref{eq:opinf_proposed} is a consequence of the orthogonality $\langle g_j,\phi_i\rangle=0$ between the basis and the projected lifting functions. Such a term would have remained if we had used the raw lifting functions $\bar{g}_j$ in the state approximation \eqref{eq:reconstruction}.
\end{remark}

Following the OpInf framework, we learn the operators in \eqref{eq:opinf_proposed} in a data-driven manner. To this end, we first compute the time series of reduced states and their time derivatives
\begin{equation}
    \{a(t_1), \dots, a(t_K)\} \subset \mathbb{R}^r, \quad \{\dot{a}(t_1), \dots, \dot{a}(t_K)\} \subset \mathbb{R}^r,
\end{equation}
corresponding to the state approximation \eqref{eq:reconstruction} applied to the full-order data snapshots. By leveraging the orthogonality between the basis and the projected lifting functions, the reduced states are simply computed as
\begin{equation}
\label{eq:projected_states_homogeneous}
    a(t_k)=\left[\left\langle y(\cdot, t_k), \phi_1\right\rangle, \dots, \left\langle y(\cdot, t_k), \phi_r\right\rangle\right]^T, \quad k=1, \ldots, K,
\end{equation}
and the corresponding derivatives $\dot{a}(t_k)$ for $k=1, \ldots, K$ are approximated by applying finite differences to the time series of reduced states, as in standard OpInf. When temporal derivatives of the full-order state are available from the PDE solver, they may be projected onto the POD basis as in~\eqref{eq:projected_states_homogeneous} to yield the reduced time derivatives.
We then infer the operators $A, H$, and $c$ so that the ROM structure derived in \eqref{eq:opinf_proposed} best fits this data in the least-squares sense, which can be formulated as the linear regression problem
\begin{equation}
\label{eq:opinf_minimization_proposed}
    \min _{A, H, c} \sum_{k=1}^K\left\|A \xi(t_k)+H(\xi(t_k) \otimes \xi(t_k))+c-\dot a(t_k)\right\|_F^2,
\end{equation}
where the augmented reduced state is given by $\xi(t_k)=[a(t_k)^T, u(t_k)^T]^T$. As in standard OpInf approaches, penalty terms involving the Frobenius norms of the inferred operators can be appended to the least-squares objective.

At inference, the reduced state at time $t_0$ is obtained from $a_i(t_0) = \langle y(\cdot, t_0), \phi_i \rangle$ for $i = 1, ..., r$, where $y(\cdot, t_0)$ is the given initial full-order state. Once the reduced state trajectory has been obtained from \eqref{eq:opinf_proposed}, the full-order state is reconstructed using \eqref{eq:reconstruction}.

We have therefore obtained a data-driven ROM that preserves the operator structure induced by boundary actuation while ensuring consistency with the prescribed boundary conditions.

\section{Experiments}
\label{sec:experiments}

\subsection{Heat Equation}
\label{sec:exp:heat}

We first evaluate the proposed framework on the one-dimensional heat equation with boundary actuation
\begin{equation}
    \frac{\partial y}{\partial t} = \frac{\partial^2 y}{\partial x^2}, 
    \quad x \in [0, 1], \quad t \in [0, 1],
\end{equation}
subject to Dirichlet BCs
\begin{equation}
\label{eq:heat_bcs}
    y(0, t) = u(t), \quad y(1, t) = u(t),
\end{equation}
where $u(t) \in \mathbb{R}$ denotes a time-varying control signal applied at both boundaries. The initial condition is chosen as
\begin{equation}
    y(x,0) =
    \left(
    e^{5(x-1)} + e^{-5 x} - e^{-5}
    \right) u(0).
\end{equation}

Assume we have collected training data in the form of a trajectory $\{y(\cdot, t_k), u(t_k)\}_{k=1}^K$. To handle the inhomogeneous BCs within a physically consistent subspace, we follow the methodology outlined in Section~\ref{sec:method} and we first select the lifting function $\bar{g}(x) = 1$, which satisfies the Dirichlet BCs \eqref{eq:heat_bcs} with $u(t)=1$. We then compute the homogeneous POD basis $\{\phi_i\}_{i=1}^r$ of the homogeneous remainder snapshots $\bar{y}_H(\cdot, t_k) = y(\cdot, t_k) - \bar{g} u(t_k)$, and we calculate the projected lifting function $g(x)$ using \eqref{eq:projected_lifting}. Finally, we approximate the full state as
\begin{equation}
\label{eq:decomposition_heat}
    y_r(x, t) = u(t) g(x) + \sum_{i=1}^r a_i(t) \phi_i(x).
\end{equation}

Since the heat equation is linear, the reduced-order dynamics \eqref{eq:opinf_proposed} of the proposed lifted OpInf model take the form
\begin{equation}
\label{eq:heat_opinf_proposed}
    \frac{d a}{d t} = A \xi(t) + c,
\end{equation}
where the reduced state is $a=[a_1, \ldots, a_r]^T \in \mathbb{R}^r$, its augmented counterpart is $\xi=[a^T, u]^T \in \mathbb{R}^{r+1}$,
and the operators $A \in \mathbb{R}^{r \times (r+1)}$ and $c \in \mathbb{R}^r$ are inferred from data by solving the least-squares problem~\eqref{eq:opinf_minimization_proposed}. In the present case of a linear PDE, the lifted OpInf reduced dynamics share the same structure as the ad hoc OpInf reduced dynamics. Nonetheless, the state approximation \eqref{eq:decomposition_heat} employed by the lifted model is guaranteed to satisfy the Dirichlet BCs, thereby maintaining consistency with the underlying boundary-actuated system.

We compare the proposed lifted OpInf model with the following methods in our numerical experiments.
For conciseness, the same notation $\{\phi_i\}_{i=1}^r$ is used for all bases. The \textit{original basis} is computed from the raw inhomogeneous snapshots $\{y(\cdot, t_k)\}_{k=1}^K$, while the \textit{lifted basis} is computed from the homogeneous remainder snapshots $\{\bar{y}_H(\cdot, t_k)\}_{k=1}^K$.

\begin{itemize}[leftmargin=*]
    \item \textbf{Projection (original)}: Projection of the full-order states onto the original POD basis $\{\phi_i\}_{i=1}^r$ obtained directly from the raw snapshots. 
    Specifically, we compute the reduced coefficients as $a_i = \langle \phi_i, y\rangle$ and reconstruct the projected state via $y_r = \sum_{i=1}^r a_i \phi_i$. 
    Although this projection does not guarantee satisfaction of the inhomogeneous BCs, it provides a bound for the best performance achievable by any Galerkin or OpInf ROM that uses this basis.

    \item \textbf{Intrusive (original)}: Galerkin ROM constructed from the original POD basis without lifting the boundary conditions. We adopt a discretize-then-project formulation in which boundary actuation enters through Galerkin projection of the spatially-discretized PDE~\cite{benner2015survey, star2021reduced}. Let $Y(t)$ denote the values of $y(\cdot, t)$ at interior grid points, and let $\Phi$ be the matrix whose columns are the discretized POD modes $\{\phi_1, \dots, \phi_r\}$ at the same points and normalized accordingly. A finite-difference discretization of the heat equation with nonzero Dirichlet BCs gives the semi-discrete system $\dot Y = A_{\text{FOM}} \ Y + B_{\text{FOM}} \ u(t)$; projecting both operators onto the POD basis then yields the ROM $\dot a = A a + B u(t)$ with $A = \Phi^\top A_{\text{FOM}} \ \Phi$, $B = \Phi^\top B_{\text{FOM}}$, and the reconstructed state is $Y_r = \Phi a$. 

    \item \textbf{OpInf (original)}: Standard OpInf ROM~\cite{qian2022reduced} constructed using the original POD basis without explicit treatment of boundary actuation. Its structure is based on the PDE with homogeneous boundary conditions, augmented with an ad hoc control input term to account for boundary actuation, as described in Section~\ref{sec:opinf}.
    The reduced dynamics are modeled as
    \begin{equation}
        \frac{d a}{d t} = A a(t) + B u(t) + c,
    \end{equation}
    with parameters $A \in \mathbb{R}^{r \times r}$, $B \in \mathbb{R}^{r \times 1}$, and $c \in \mathbb{R}^r$ identified from data. The state is reconstructed as $y_r = \sum_{i=1}^r a_i \phi_i$.

    \item \textbf{Projection (lifted)}: Projection of the full-order states using the proposed lifted decomposition \eqref{eq:decomposition_heat} and the POD basis $\{\phi_i\}_{i=1}^r$ computed from remainder states. 
    Specifically, we compute the reduced coefficients as $a_i = \langle \phi_i, y\rangle$ and reconstruct the projected state via $y_r = u g + \sum_{i=1}^r a_i \phi_i$. 
    This projection guarantees boundary condition satisfaction, as shown in Theorem~\ref{thm:bc_satisfaction}, and provides a benchmark for best attainable performance by any lifted Galerkin or OpInf ROM using such basis.

    \item \textbf{Intrusive (lifted)}: Lifted Galerkin ROM~\cite{hakenberg2013model} derived from the lifted decomposition \eqref{eq:decomposition_heat}, yielding reduced dynamics of the form~\eqref{eq:heat_opinf_proposed} without the offset term $c$. The operator  $A$ is calculated explicitly from the residual constraint \eqref{eq:residual_condition_inhomogeneous} with $\mathcal{L}$ the Laplacian operator and $\mathcal{Q} = 0$. The terms $\mathcal{L}\phi_j$ and $\mathcal{L}g$ are approximated using finite differences that incorporate the zero endpoint values of $\phi_j$ and nonzero endpoint values of $g$, and the inner products are evaluated using trapezoidal quadrature.

    \item \textbf{OpInf (lifted)}: The proposed OpInf ROM based on the physically consistent lifted formulation of Section~\ref{sec:method}. The model uses the lifted decomposition \eqref{eq:decomposition_heat} with the reduced dynamics~\eqref{eq:heat_opinf_proposed}.
\end{itemize}

To evaluate these methods, we consider the following time-varying boundary control functions:
\begin{subequations} \label{eq:heat_control_functions}
\begin{align}
    u_1(t) &= 1 + t(1 - t), \\
    u_2(t) &= 1 + \tfrac{1}{4} \sin(4\pi t), \\
    u_3(t) &= 1 + 5 \sin(10\pi t) e^{-2t}, \\
    u_4(t) &= 1 + 10 \sin(50\pi t^2), \\
    u_5(t) &= 1 - \tfrac{1}{2} \sin(3\pi t).
\end{align}
\end{subequations}
High-fidelity trajectories corresponding to these control functions are generated by discretizing the spatial domain with a Chebyshev collocation method using $N+1=65$  Chebyshev–Lobatto spatial nodes and integrating the resulting system with the adaptive RK45 method implemented in SciPy's \texttt{solve\_ivp}. Snapshots are interpolated onto a common uniform spatial grid and recorded at intervals of $\Delta t = 10^{-3}$, yielding five trajectories with $K=1001$ snapshots each.

We use the official Operator Inference Python package\footnote{\url{https://github.com/Willcox-Research-Group/rom-operator-inference-Python3}} to implement the original OpInf and original intrusive Galerkin ROM baselines. 
The least-squares problem~\eqref{eq:opinf_minimization_proposed} for the lifted OpInf formulation is solved by forming the normal equations and solving the resulting linear systems using \texttt{numpy.linalg.solve}. Both original and lifted OpInf formulations are solved without regularization, and the time derivative of the reduced state in their regression objectives is obtained using 6th-order finite differences. During inference, the ROMs are integrated using the variable-order backward differentiation formula (BDF) method implemented in SciPy's \texttt{solve\_ivp}.

For each method, we construct five separate models, with each model trained on the trajectory generated by one boundary control function in \eqref{eq:heat_control_functions}. Each trained model is then evaluated on all five boundary control functions. This yields 5 evaluations in which the boundary control function used for evaluation is the same as that used for training, and 20 evaluations in which it is different. We report the training error as the average over the 5 same-control evaluations and the testing error as the average over the 20 different-control evaluations.

Table~\ref{tab:heat_solver_failure} first reports the solver failure rates for the original and lifted OpInf formulations over the evaluated training and testing cases. Both formulations exhibit increasing failure rates as the number of POD modes grows, but the lifted OpInf formulation consistently achieves lower failure rates whenever failures occur. In particular, the lifted OpInf model can be successfully simulated for every case in which the original OpInf model succeeds. This suggests that the lifted formulation improves numerical robustness, even though the least-squares regression is overdetermined for both formulations\footnote{In this heat equation example, the number of equations contributing to the square objective is $rK$ while the number of unknowns is $r^2 + 2r$ and $r(r+1) + 2r$ for the original and lifted formulations, respectively. Since $K \gg r$, the least-squares regression is vastly overdetermined.}. To ensure a fair comparison of approximation accuracy, all subsequent error results are computed over the same subset of cases for which the original OpInf model can be successfully simulated, with the lifted OpInf model also yielding finite predictions in all these cases.

\begin{figure}[t]
    \centering
    
    \begin{subfigure}{0.49\linewidth}
        \centering
        \includegraphics[width=\linewidth]
        {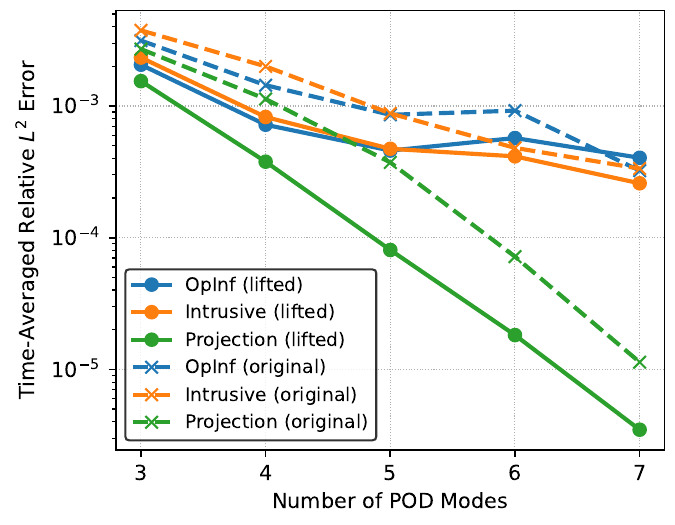}
        \caption{Training error.}
        \label{fig:heat_train_error_vs_pod_modes}
    \end{subfigure}
    \hfill
    \begin{subfigure}{0.49\linewidth}
        \centering
        \includegraphics[width=\linewidth]
        {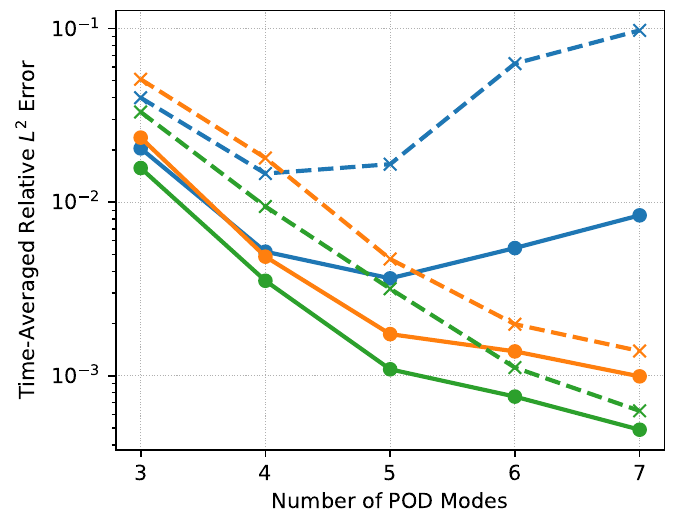}
        \caption{Testing error.}
        \label{fig:heat_test_error_vs_pod_modes}
    \end{subfigure}
    
    \caption{Time-averaged relative $L^2$ prediction error with respect
    to the number of POD modes for the heat equation.}
    \label{fig:heat_error_vs_pod_modes}
\end{figure}

\begin{figure}[t]
    \centering
    \includegraphics[width=0.5\linewidth]{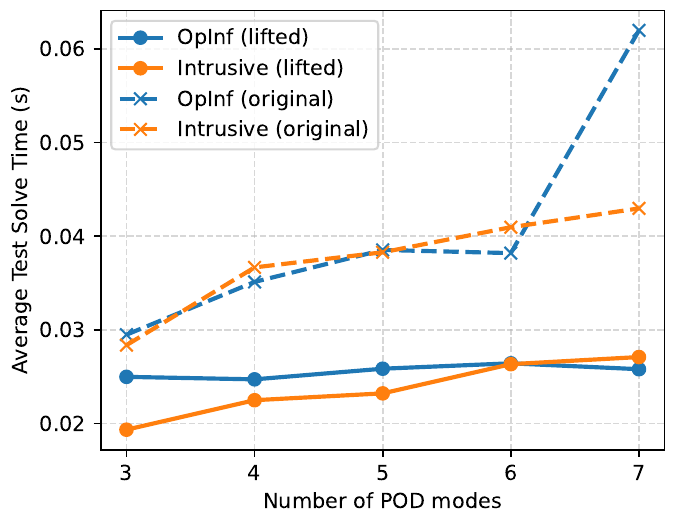}
    \caption{Test-set solver time for the heat equation with different numbers of POD modes.}
    \label{fig:heat_solver_time_vs_pod_modes}
\end{figure}

\begin{figure}[t]
    \centering
    
    \begin{subfigure}{0.49\linewidth}
        \centering
        \includegraphics[width=\linewidth]{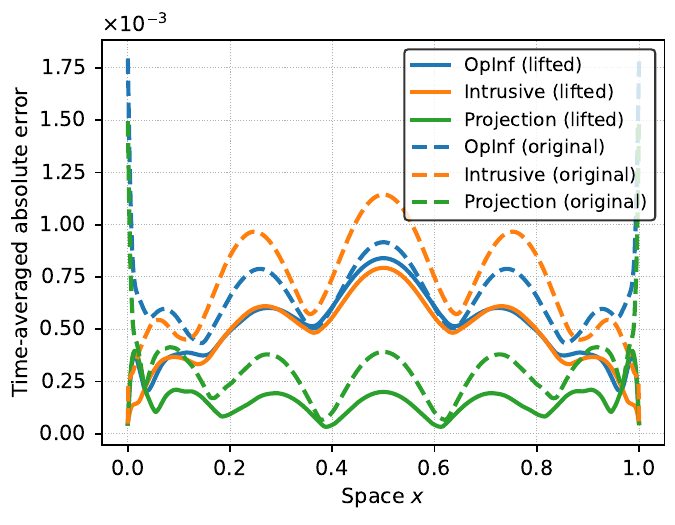}
        \caption{Training error.}
        \label{fig:heat_train_state_error}
    \end{subfigure}
    \hfill
    \begin{subfigure}{0.49\linewidth}
        \centering
        \includegraphics[width=\linewidth]{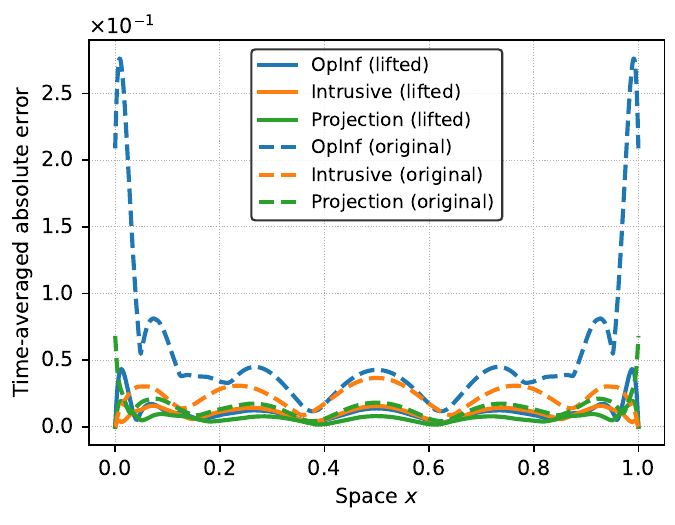}
        \caption{Testing error.}
        \label{fig:heat_test_state_error}
    \end{subfigure}
    
    \caption{Time-averaged absolute prediction error versus spatial location for the heat equation, averaged over all POD dimensions.}
    \label{fig:heat_state_error}
\end{figure}

Fig.~\ref{fig:heat_error_vs_pod_modes} displays the relative $L^2$ prediction error for different POD subspace dimension $r$ from 3 to 8. The proposed lifting scheme consistently reduces the training and testing errors of the projection, OpInf, and intrusive ROM approaches. Fig.~\ref{fig:heat_solver_time_vs_pod_modes} reports the average test-set solve time as the number of POD modes increases. The lifted OpInf and intrusive ROMs both exhibit relatively stable solve times across the tested POD subspace dimensions, while the corresponding original formulations show larger increases as the subspace dimension grows. These results imply that, for this heat equation problem, the proposed lifting scheme can improve online prediction efficiency in addition to accuracy and robustness. Finally, the point-wise $L^1$ prediction error over the spatial domain is visualized in Fig.~\ref{fig:heat_state_error}, averaged over all tested POD subspace dimensions. Thanks to strict boundary-condition satisfaction, the lifted projection, intrusive ROM, and OpInf methods substantially reduce the error near the boundary.

\begin{table}[t]
\centering
\caption{Solver failure rate (\%) on train and test data for the heat equation at different numbers of POD modes.}
\label{tab:heat_solver_failure}
\begin{tabular}{lcccccc}
\toprule
Method / POD Modes & 3 & 4 & 5 & 6 & 7 & 8 \\
\midrule
OpInf (lifted)   & 0.00 & 0.00 & 0.00 & 20.00 & 40.00 & 80.00 \\
OpInf (original) & 0.00 & 0.00 & 0.00 & 40.00 & 60.00 & 100.00 \\
\bottomrule
\end{tabular}
\end{table}

\subsection{Burgers' Equation}
\label{sec:exp:burgers}

We then evaluate the proposed operator inference framework on the viscous Burgers’ equation with boundary control and interior forcing
\begin{equation}
    \frac{\partial y}{\partial t} + y\,\frac{\partial y}{\partial x}
    = \nu\,\frac{\partial^2 y}{\partial x^2} + s(x)\,w_3(t),
    \quad x \in [0, 1], \quad t \in [0, 2],
\end{equation}
subject to Dirichlet boundary conditions
\begin{equation}
\label{eq:burgers_bcs}
    y(0, t) = w_1(t), \quad y(1, t) = w_2(t),
\end{equation}
where $w_1(t), w_2(t), w_3(t) \in \mathbb{R}$ are time-dependent control signals, and the interior forcing is specified through the spatial profile 
$s(x) = \cosh((x-0.5)/0.05)^{-1}$. The initial condition is set consistently with the boundary values as
\begin{equation}
    y(x,0)
    =
    w_2(0)
    +
    \frac{1}{2}\bigl(w_1(0)-w_2(0)\bigr)
    \left[
    1-\tanh\!\left(\frac{x-0.3}{0.1}\right)
    \right].
\end{equation}

Assume we have collected training data in the form of a trajectory $\{y(\cdot, t_k), w_1(t_k), w_2(t_k), w_3(t_k)\}_{k=1}^K$. To handle the inhomogeneous BCs within a physically consistent subspace, we first select the lifting functions $\bar{g}_1(x) = 1-x$ and $\bar{g}_2(x) = x$, which satisfy the Dirichlet BCs \eqref{eq:burgers_bcs} for $[w_1, w_2]=[1,0]$ and $[0,1]$, respectively. We then compute the homogeneous POD basis $\{\phi_i\}_{i=1}^r$ of the homogeneous remainder snapshots $\bar{y}_H(\cdot, t_k) = y(\cdot, t_k) - \bar{g}_1 w_1(t_k) - \bar{g}_2 w_2(t_k)$, and we calculate the projected lifting functions $g_1(x)$ and $g_2(x)$ using \eqref{eq:projected_lifting}. Finally, we approximate the full state as
\begin{equation}
\label{eq:burgers}
    y_r(x, t) = w_1(t) \, g_1(x) + w_2(t) \, g_2(x) + \sum_{i=1}^r a_i(t) \phi_i(x),
\end{equation}

The reduced-order dynamics of the proposed lifted OpInf model then take the form
\begin{equation}
\label{eq:burgers_opinf_proposed}
    \frac{da}{dt} = A \xi(t) + H(\xi(t) \otimes \xi(t)) + b w_3(t) + c,
\end{equation}
where the reduced state is $a=[a_1, \ldots, a_r]^T \in \mathbb{R}^r$, its augmented counterpart is $\xi=[a^T, w_1, w_2]^T \in \mathbb{R}^{r+2}$,
the operators have shapes $A \in \mathbb{R}^{r \times (r+2)}$, $H \in \mathbb{R}^{r \times (r+2)(r+2)}$, and $b, c \in \mathbb{R}^r$, and the quadratic contraction is defined component-wise by
\begin{equation}
    [H(\xi \otimes \xi)]_i
    = \sum_{j=1}^{r+2}\sum_{k=1}^{r+2} H_{i,(j-1)(r+2)+ k}\,\xi_j\,\xi_k,
    \quad i = 1,\dots,r.
\label{eq:quadratic_contraction}
\end{equation}
Finally, the corresponding least-squares regression problem for the operators is given by
\begin{equation}
\label{eq:burgers_proposed_opt}
    \min_{A, H, b, c} \sum_{k=1}^K\left\|A \xi(t_k) + H(\xi(t_k) \otimes \xi(t_k))
 + b w_3(t_k) + c-\dot{a}(t_k)\right\|_F^2.
\end{equation}
\begin{remark}
Unlike boundary control, the interior forcing term $s(x) w_3(t)$ does not change the admissible solution space of the PDE, so $w_3$ need not enter the lifting decomposition in \eqref{eq:burgers}. Galerkin projection of this term instead yields a standard additive component linear in $w_3$ in the reduced dynamics \eqref{eq:burgers_opinf_proposed}, whereas the boundary inputs $w_1$ and $w_2$ include linear, bilinear, and quadratic contributions through $\xi$. 
\end{remark}

We compare the lifted method with the original OpInf~\cite{kramer2024learning} based on the POD approximation in \eqref{eq:pod_expansion}, with the ROM dynamics taking the form
\begin{equation}
    \frac{da}{dt} = A a(t) + H(a(t) \otimes a(t))
    + b_1 w_1(t) + b_2 w_2(t) + b_3 w_3(t) + c,
\end{equation}
where the operators have shapes $A \in \mathbb{R}^{r \times r}$, $H \in \mathbb{R}^{r \times r^2}$, and $b_1, b_2, b_3, c \in \mathbb{R}^r$. The regression problem is given by
\begin{equation}
\label{eq:burgers_opinf_opt}
    \min_{A, H, b_1, b_2, b_3, c} \sum_{k=1}^K\left\|A a(t_k) + H(a(t_k) \otimes a(t_k))
    + b_1 w_1(t_k) + b_2 w_2(t_k) + b_3 w_3(t_k) + c-\dot{a}(t_k)\right\|_F^2.
\end{equation}

It was observed in Section~\ref{sec:exp:heat} that the lifting approach enhances the well-posedness and robustness of the inferred OpInf model. Since ridge regression is a common approach to address these issues, here we apply this regularization to the least-squares regression problems for both the lifted and original OpInf formulations. Specifically, letting $B=b$ for the lifted OpInf and $B=[b_1,b_2,b_3]$ for the original OpInf, we add the following Frobenius-norm penalty terms to~\eqref{eq:burgers_proposed_opt} and~\eqref{eq:burgers_opinf_opt},
\begin{equation}
\label{eq:burgers_regularization}
\mathcal{R}(A, H, B, c)= \begin{cases}
\|H\|_F^2+\|A\|_F^2+\|B\|_F^2+\|c\|_2^2, & \text { if full regularization (OpInf),} \\
\|H\|_F^2, & \text { if quadratic-only regularization (PIR-OpInf),}
\end{cases}
\end{equation}
multiplied by a scalar ridge regression weight $\lambda$. The quadratic-only regularization is inspired by prior literature on physics-informed regularization for quadratic and cubic OpInf models (PIR-OpInf) in~\cite{sawant2023physics}, which showed that solely penalizing higher-order terms can potentially enhance model stability and robustness. 

To evaluate these methods on varied trajectories, we define each input signal as a sum of multiple sine waves with varying magnitudes, frequencies, phases, and biases as follows
\begin{equation}
    w_k(t) = b_k + \sum_{m=1}^{3} a_{k,m}\,\sin\!\big(2\pi f_{k,m} t + \varphi_{k,m}\big),
    \qquad k \in \{1,2,3\},
\end{equation}
where the parameters corresponding to the boundary control inputs $w_1$ and $w_2$ are sampled from
\begin{equation}
    a_{k,m} \sim \mathcal{U}[0.1, 1.0], \quad
    f_{k,m} \sim \mathcal{U}[0.2, 5.0], \quad
    \varphi_{k,m} \sim \mathcal{U}[0, 2\pi], \quad
    b_k \sim \mathcal{U}[-0.5, 0.5], 
\end{equation}
for $k \in \{1,2\}$, while those corresponding to the interior forcing term $w_3$ are sampled from
\begin{equation}
    a_{3,m} \sim \mathcal{U}[0.1, 3.0], \quad
    f_{3,m} \sim \mathcal{U}[0.2, 1.0], \quad
    \varphi_{3,m} \sim \mathcal{U}[0, 2\pi], \quad
    b_3 \sim \mathcal{U}[-1.0, 1.0].
\end{equation}
High-fidelity trajectories corresponding to these input functions are generated by discretizing the spatial domain with a Chebyshev collocation method using $N+1=65$ Chebyshev–Lobatto spatial nodes and integrating the resulting system with the adaptive RK45 method implemented in SciPy's \texttt{solve\_ivp}. Snapshots are recorded at intervals of $\Delta t = 0.02$, yielding trajectories with $K=101$ snapshots each. 
Since the snapshots are represented on a nonuniform Chebyshev grid, POD projections, modal coefficients, and error metrics are computed using a weighted discrete inner product with Clenshaw--Curtis quadrature weights scaled to the spatial domain to approximate the standard $L^2$ inner product.
The ridge-regularized OpInf regression problems~\eqref{eq:burgers_proposed_opt} and~\eqref{eq:burgers_opinf_opt} are solved by forming the normal equations and solving the resulting dense linear systems using \texttt{numpy.linalg.solve}. The time derivative of the reduced state in the regression objectives is obtained using 4th-order finite differences. 
During inference, the ROMs are integrated with SciPy's \texttt{solve\_ivp} using BDF, with RK45 used as a fallback if BDF fails.

We generate 100 trajectories each for training, validation, and testing. For both the lifted and original OpInf methods, we train models using 20 ridge regression weights $\lambda$ sampled log-uniformly in $[10^{-4},10^2]$. The best regularization parameter $\lambda^*$ is selected separately for each method based on its prediction performance on the validation set, and the corresponding model is then evaluated on the test set.

\begin{figure}[t]
    \centering
    
    \begin{subfigure}{0.49\linewidth}
        \centering
        \includegraphics[width=\linewidth]{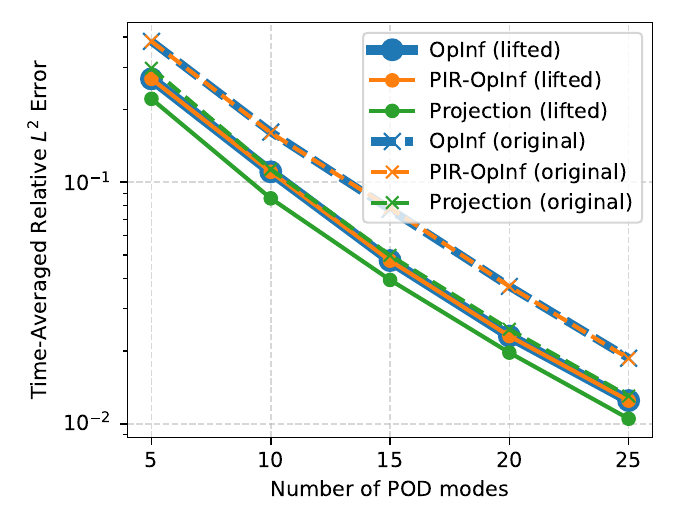}
        \caption{Time-averaged $L^2$ error.}
    \end{subfigure}
    \hfill
    \begin{subfigure}{0.49\linewidth}
        \centering
        \includegraphics[width=\linewidth]{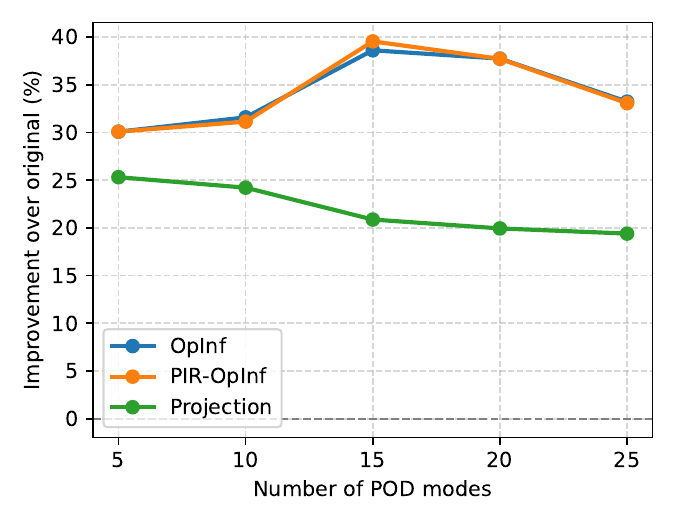}
        \caption{Percentage error reduction from the original to the lifted formulation.}
    \end{subfigure}
    
    \caption{Time-averaged $L^2$ prediction error with respect to number of POD modes for Burgers' equations.}
    \label{fig:burgers_error_pod}
\end{figure}

Fig.~\ref{fig:burgers_error_pod} reports the $L^2$ prediction error of all methods for different numbers of POD modes, together with the error reduction achieved by the lifted formulation relative to the original one. The projection methods are similar to those described in Section \ref{sec:exp:heat}, adapted to the decomposition \eqref{eq:burgers}. The lifted OpInf method consistently outperforms the original OpInf formulation across different POD subspace dimensions, using either full or quadratic-only regularization, which yield comparable results. In addition, the projection error curves show that the proposed lifting scheme yields a lower best-attainable ROM error for a given number of POD modes. This improvement is likely due to the consistency of the lifted formulation with the physical boundary conditions. 

\begin{figure}[t!]
    \centering
    \includegraphics[width=0.5\linewidth]{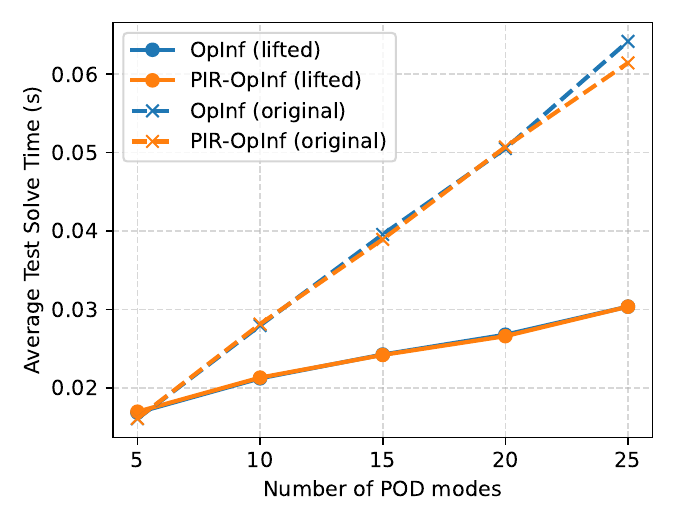}
    \caption{Test-trajectory solver time for the Burgers' equation with different numbers of POD modes.}
    \label{fig:burgers_time_modes}
\end{figure}

Next, Table~\ref{tab:burgers_solver_failure} reports the number of solver failures for the lifted and original OpInf methods during validation across all 20 tested regularization parameters $\lambda \in [10^{-4},10^2]$. The results show that the lifted method significantly improves solver stability irrespective of different regularization, especially for larger numbers of POD modes, where the original methods begin to exhibit validation failures. Fig.~\ref{fig:burgers_time_modes} compares the solver prediction time of the lifted and original formulations on the test set. The lifted formulations exhibit substantially more favorable scaling with the POD subspace dimension, becoming faster than the original formulations for more than 10 POD modes.

\begin{table}[t]
\centering
\caption{Solver failure rate (\%) on validation data for Burgers' equations at different numbers of POD modes.}
\label{tab:burgers_solver_failure}
\begin{tabular}{lccccc}
\toprule
Method / POD Modes & 5 & 10 & 15 & 20 & 25 \\
\midrule
OpInf (lifted)     & 0.00 & 0.00 & 0.00 & 0.00 & 0.00 \\
PIR-OpInf (lifted) & 0.00 & 0.00 & 0.00 & 0.00 & 0.00 \\
OpInf (original)     & 0.00 & 0.00 & 0.00 & 0.00 & 0.45 \\
PIR-OpInf (original) & 0.00 & 0.00 & 0.00 & 2.35 & 30.15 \\
\bottomrule
\end{tabular}
\end{table}

\begin{figure}[t!]
    \centering
    \includegraphics[width=0.95\linewidth]{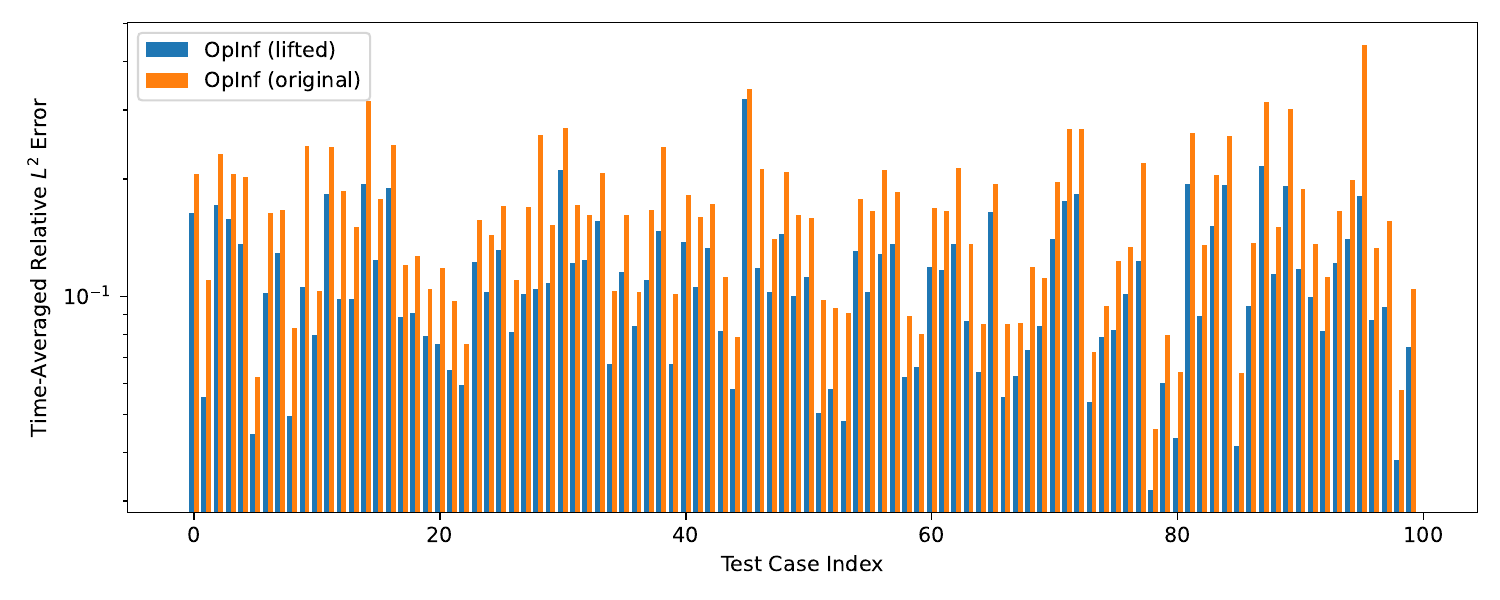}
    \caption{Per-case test-set relative $L^2$ prediction error for Burgers' equations with 10 POD modes.}
    \label{fig:burger_error_stats}
\end{figure}

\begin{figure}[t!]
    \centering
    \includegraphics[width=0.95\linewidth]{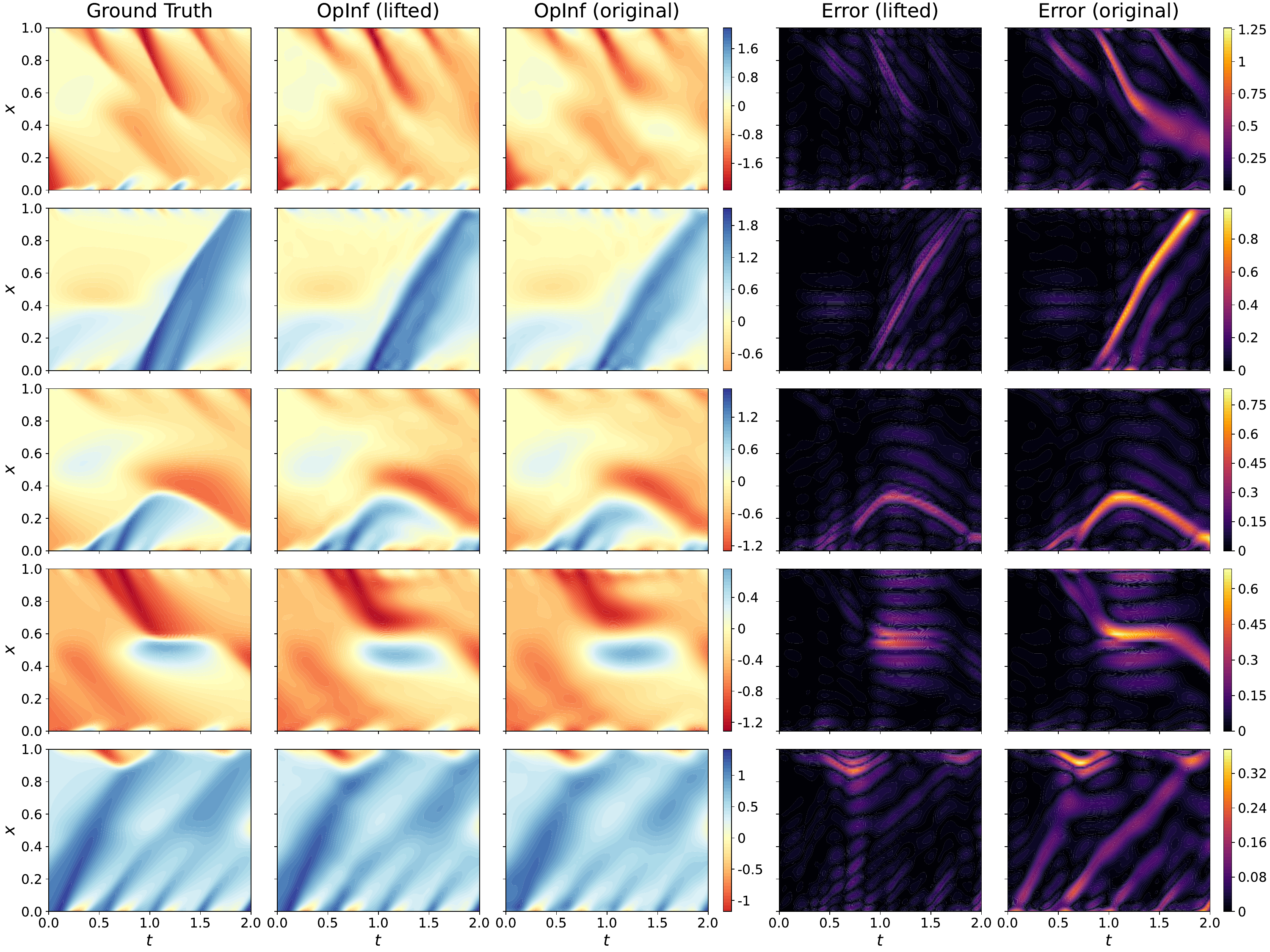}
    \caption{Visualization of representative test trajectories and pointwise absolute errors for Burgers' equation with 10 POD modes. The first three columns show the reference, lifted OpInf, and original OpInf solution trajectories. The final two columns show the corresponding pointwise absolute solution errors.}
    \label{fig:burger_visualization}
\end{figure}

Finally, Fig.~\ref{fig:burger_error_stats} shows the per-trajectory test-set $L^2$ prediction error for the lifted and original OpInf methods, and Fig.~\ref{fig:burger_visualization} visualizes representative test trajectories, under the setting of 10 POD modes with full regularization. It can be observed that the lifted method consistently reduces the prediction error compared to the original OpInf formulation across all testing instances. These results further demonstrate the effectiveness of the proposed lifting scheme and its ability to handle diverse Burgers' equation settings with both boundary control and source terms.

\subsection{2D Lid-Driven Cavity}
\label{sec:exp:cavity}

As a last example, we consider the canonical problem of incompressible flow in a lid-driven cavity, where the lid velocity is prescribed as an arbitrary time-varying input. The flow takes place in the two-dimensional unit square domain $\Omega = [0,1]^2$ and is governed by the incompressible Navier-Stokes equations, written in non-dimensional form as
\begin{subequations} \label{eq:NS}
\begin{align}
\frac{\partial \mathbf{v}}{\partial t} + (\mathbf{v} \cdot \nabla)\mathbf{v} &= -\nabla p + \frac{1}{Re} \Delta \mathbf{v}, \quad \mathbf{x} \in \Omega, \label{eq:NS1} \\
\nabla \cdot \mathbf{v} &= 0, \quad \mathbf{x} \in \Omega, \label{eq:NS2}
\end{align}
\end{subequations}
where $\mathbf{x} = (x,y)$ is the spatial coordinate, $\mathbf{v} = (v,w)$ denotes the velocity, $p$ is the pressure, and $Re$ is the Reynolds number selected herein to be 100. No-slip boundary conditions are imposed on all walls. The horizontal lid is prescribed a uniform tangential velocity $u(t) \in \mathbb{R}$ while the remaining boundaries are stationary, resulting in the time-dependent Dirichlet BCs
\begin{subequations} \label{eq:NS_BCs}
\begin{align}
\mathbf{v}(x,0) = (0,0), \quad \mathbf{v}(x,1) &= (u(t),0), \quad x \in [0,1], \\
\mathbf{v}(0,y) = (0,0), \quad \mathbf{v}(1,y) &= (0,0), \quad y \in [0,1].
\end{align}
\end{subequations}

The Navier-Stokes equations are discretized in the open-source solver OpenFOAM \cite{openfoam_opencfd} using a graded $100\times100$ finite-volume mesh, with increased resolution near the cavity walls, and integrated using the transient solver \texttt{icoFoam} with a time step of $\Delta t=0.002$. An approximately steady solution $\mathbf{v}_s$ corresponding to constant lid velocity $u=1$ is first computed and subsequently used as the initial condition for the time-dependent simulation. In the latter, the lid velocity is prescribed as
\begin{equation}
u(t)=\sum_{m=1}^{5} a_m \sin\left(\omega_0 f_m t+\varphi_m\right),
\qquad \omega_0=\frac{2\pi}{T},
\end{equation}
where the initial amplitudes, frequency factors, and phases are sampled as
\begin{equation}
a_m \sim \mathcal{U}[0,1], \quad
f_m \sim \mathcal{U}[1,20], \quad
\varphi_m \sim \mathcal{U}[0,2\pi],
\qquad m=1,\ldots,5.
\end{equation}
The amplitudes are then jointly rescaled so that $u(0)=1$, ensuring consistency between the steady initial condition and the time-dependent boundary condition. The simulation is performed over $t\in[0,T]$ with $T=200$, and the velocity and lid velocity are recorded every 0.1 time units, yielding a trajectory $\{\mathbf{v}(\cdot, t_k), u(t_k)\}_{k=1}^K$ with $K = 1001$ snapshots.

As in \cite{benner2020operator2}, we first reformulate the Navier-Stokes equations \eqref{eq:NS} into a form that does not explicitly depend on the pressure, so that OpInf models can be constructed from velocity snapshots alone. Taking the divergence of the momentum equation \eqref{eq:NS1} and applying the incompressibility condition \eqref{eq:NS2} leads to the pressure Poisson equation
\begin{equation}
\Delta p = - \nabla \cdot [(\mathbf{v} \cdot \nabla)\mathbf{v}],
\end{equation}
where we have used the commutativity of the divergence operator with the time derivative and the Laplace operator. Inserting the corresponding solution for the pressure back into \eqref{eq:NS1}, the Navier-Stokes equations are reduced into a quadratic form that only involves the velocity field,
\begin{equation}
\frac{\partial \mathbf{v}}{\partial t} + \mathbb{P}[ (\mathbf{v} \cdot \nabla)\mathbf{v} ] = \frac{1}{Re} \Delta \mathbf{v}, \quad \mathbf{x} \in \Omega, \label{eq:vel}
\end{equation}
where $\mathbb{P}$ is a linear operator known as the Leray projector and formally defined as $\mathbb{P}(\mathbf{s}) = \mathbf{s} - \nabla \Delta^{-1} (\nabla \cdot \mathbf{s})$. The boundary conditions associated with \eqref{eq:vel} remain the same as those in \eqref{eq:NS_BCs}. 

Following our proposed methodology, we handle the inhomogeneous BCs within a physically consistent subspace by lifting the solution. We choose as lifting function $\bar{\mathbf{g}}$ the steady-state velocity field $\mathbf{v}_s$ satisfying the inhomogeneous Dirichlet BCs \eqref{eq:NS_BCs} with $u = 1$. We then compute the homogeneous POD basis $\{\boldsymbol{\phi}_i\}_{i=1}^r$ of the homogeneous remainder snapshots $\bar{\mathbf{v}}_H(\cdot, t_k) = \mathbf{v}(\cdot, t_k) - \bar{\mathbf{g}} \, u(t_k)$, and we calculate the projected lifting function $\mathbf{g}$ using \eqref{eq:projected_lifting}. Note that $\mathbf{g}$ is divergence-free because both $\bar{\mathbf{g}}$ and the POD modes are divergence-free. Finally, we approximate the full state as
\begin{equation}
\label{eq:NS_decomposition}
\mathbf{v}_r(\mathbf{x},t) = \mathbf{g}(\mathbf{x}) u(t) + \sum_{i=1}^r a_i(t)\boldsymbol{\phi}_i(\mathbf{x}).
\end{equation}

Due to the quadratic structure of the velocity-only equation \eqref{eq:vel} as well as the divergence-free nature of \eqref{eq:NS_decomposition}, the reduced-order dynamics of the proposed lifted OpInf model take the form
\begin{equation}
    \frac{da}{dt} = A \xi(t) + H(\xi(t) \otimes \xi(t)) + c,
\end{equation}
where the reduced state is $a=[a_1, \ldots, a_r]^T \in \mathbb{R}^r$, its augmented counterpart is $\xi=[a^T, u]^T \in \mathbb{R}^{r+1}$,
the operators have shapes $A \in \mathbb{R}^{r \times (r+1)}$, $H \in \mathbb{R}^{r \times (r+1)(r+1)}$, $c \in \mathbb{R}^r$, and the quadratic contraction is defined as in \eqref{eq:quadratic_contraction}. Finally, the operators are inferred by solving the least-squares regression problem
\begin{equation}
\label{eq:NS_proposed_opt}
    \min _{A, H, c} \sum_{k=1}^K\left\|A \xi(t_k) + H(\xi(t_k) \otimes \xi(t_k))
    + c-\dot{a}(t_k)\right\|_F^2. 
\end{equation}

We compare the lifted method with the original OpInf~\cite{benner2020operator2,kramer2024learning} based on the POD approximation in \eqref{eq:pod_expansion}, with the ROM dynamics taking the form
\begin{equation}
    \frac{da}{dt} = A a(t) + H(a(t) \otimes a(t))
    + b u(t) + c,
\end{equation}
where the operators have shapes $A \in \mathbb{R}^{r \times r}$, $H \in \mathbb{R}^{r \times r^2}$, and $b, c \in \mathbb{R}^r$. The regression problem is given by
\begin{equation}
\label{eq:NS_opinf_opt}
    \min _{A, H, b, c} \sum_{k=1}^K\left\|A a(t_k) + H(a(t_k) \otimes a(t_k))
    + b u(t_k) + c-\dot{a}(t_k)\right\|_F^2.
\end{equation}

As in the Burgers' example, we consider both full and quadratic-only regularization by augmenting the OpInf regression problems in \eqref{eq:NS_proposed_opt} and \eqref{eq:NS_opinf_opt} with the penalty terms defined in \eqref{eq:burgers_regularization}, taking $B$ to be absent for the lifted model and $B=b$ for the original one. 
The resulting ridge-regression problems are solved by forming the normal equations and solving the resulting dense linear systems using \texttt{numpy.linalg.solve}. The time derivative of the reduced state in the regression objectives is obtained using 4th-order finite differences. During inference, the ROMs are integrated with SciPy's \texttt{solve\_ivp} using BDF, with RK45 used as a fallback if BDF fails. 

We train the lifted and original OpInf models on the time interval $t\in[0,100]$, use the validation interval $t\in[100,120]$ to select the regularization parameter $\lambda^*$ from 20 logarithmically spaced candidates in $[10^{-10},10^1]$, and evaluate the resulting models on the test interval $t\in[120,200]$.

\begin{figure}[t]
    \centering
    
    \begin{subfigure}{0.49\linewidth}
        \centering
        \includegraphics[width=\linewidth]{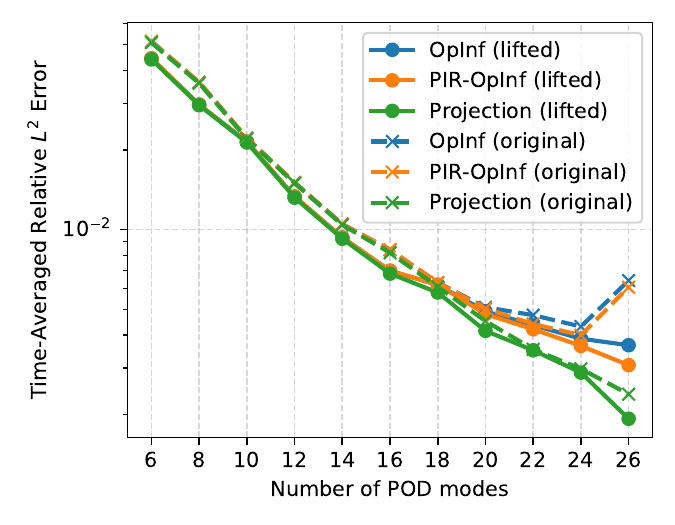}
        \caption{Time-averaged relative $L^2$ error.}
        \label{fig:cavity_error_pod_abs}
    \end{subfigure}
    \hfill
    \begin{subfigure}{0.49\linewidth}
        \centering
        \includegraphics[width=\linewidth]{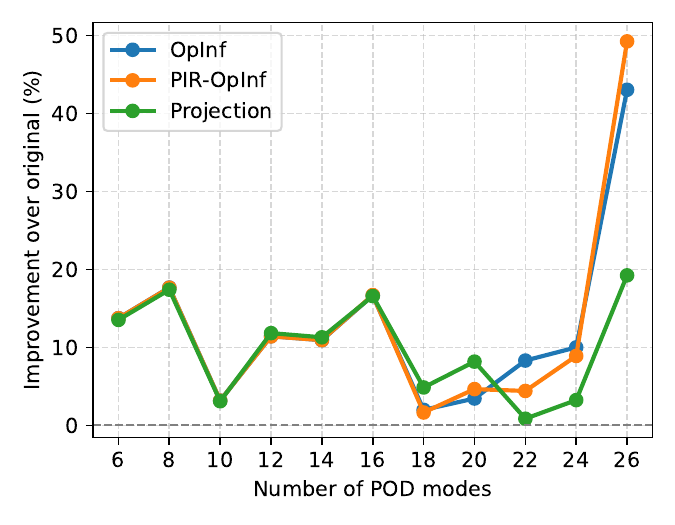}
        \caption{Percentage error reduction from the original to the lifted formulation.}
        \label{fig:cavity_error_pod_reduction}
    \end{subfigure}
    
    \caption{Time-averaged relative $L^2$ prediction error with respect to number of POD modes in the lid-driven cavity problem.}
    \label{fig:cavity_error_pod}
\end{figure}

Fig.~\ref{fig:cavity_error_pod} reports the relative $L^2$ prediction error of all methods for different numbers of POD modes, together with the error reduction achieved by the lifted formulation relative to the original one. The projection methods are similar to those described in Section \ref{sec:exp:heat}, adapted to the decomposition \eqref{eq:NS_decomposition}. The lifted method improves upon the original OpInf formulation under all settings for both full and quadratic-only regularization. Notably, the original OpInf model achieves its best performance at 24 POD modes and deteriorates at 26 modes, whereas the error of the lifted formulation continues to decrease. The projection error is also consistently reduced, indicating that the proposed lifting scheme yields a more accurate reduced representation. 

\begin{figure}[t!]
    \centering
    \includegraphics[width=0.5\linewidth]{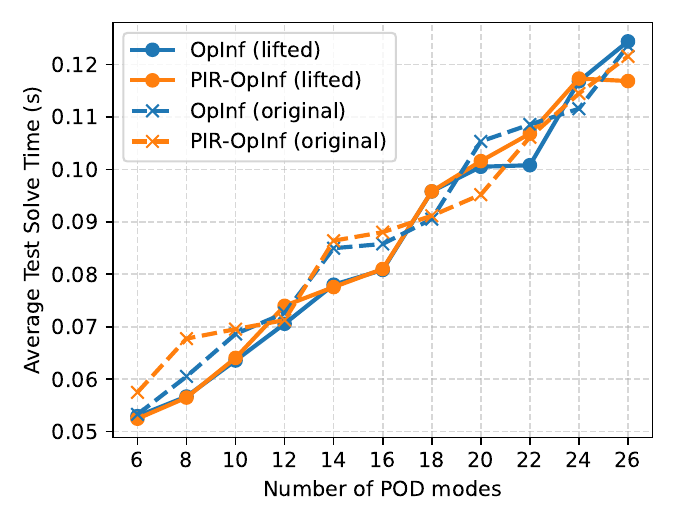}
    \caption{Testing solver time for the lid-driven cavity problem with different numbers of POD modes.}
    \label{fig:cavity_test_time_vs_pod_modes}
\end{figure}

As shown in Table~\ref{tab:cavity_solver_failure}, all formulations remain stable for low-dimensional ROMs, with no validation failures up to 12 POD modes. At higher dimensions, validation failures become more frequent, especially for the original PIR-OpInf model, while the lifted formulation generally maintains comparable or lower failure rates. Thus, the lifting scheme improves accuracy without introducing a noticeable robustness penalty. Fig.~\ref{fig:cavity_test_time_vs_pod_modes} compares the average solver prediction time for different numbers of POD modes on the test set. In this lid-driven cavity example, the lifted and original formulations have comparable solve times, and the lifted formulation does not provide a measurable improvement in online efficiency. This may be due to the additional terms in the ROM dynamics induced by the lifting transformation, including cross terms between the boundary input and the modal coefficients, as well as quadratic input terms.

\begin{table}[t]
\centering
\caption{Validation failure rate (\%) for the cavity problem at different numbers of POD modes.}
\label{tab:cavity_solver_failure}
\begin{tabular}{lccccccccccc}
\toprule
Method / POD Modes & 6 & 8 & 10 & 12 & 14 & 16 & 18 & 20 & 22 & 24 & 26 \\
\midrule
OpInf (lifted)       & 0.00 & 0.00 & 0.00 & 0.00 & 5.00 & 0.00 & 5.00 & 15.00 & 15.00 & 15.00 & 15.00 \\
PIR-OpInf (lifted)   & 0.00 & 0.00 & 0.00 & 0.00 & 5.00 & 0.00 & 5.00 & 15.00 & 15.00 & 15.00 & 20.00 \\
OpInf (original)     & 0.00 & 0.00 & 0.00 & 0.00 & 0.00 & 0.00 & 5.00 & 10.00 & 10.00 & 15.00 & 15.00 \\
PIR-OpInf (original) & 0.00 & 0.00 & 0.00 & 0.00 & 0.00 & 0.00 & 10.00 & 20.00 & 20.00 & 35.00 & 15.00 \\
\bottomrule
\end{tabular}
\end{table}

\begin{figure}[t!]
    \centering
    \includegraphics[width=0.95\linewidth]{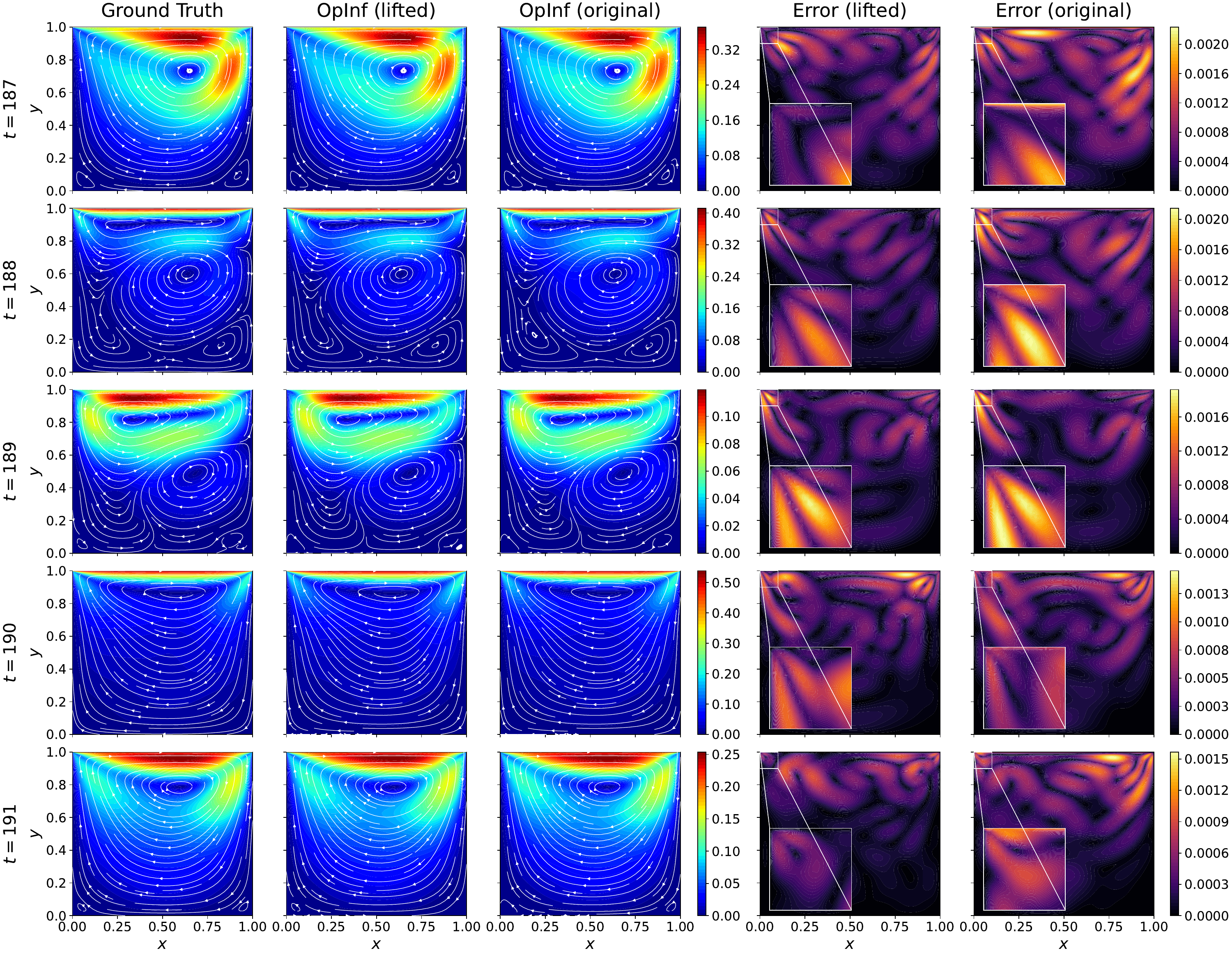}
    \caption{Representative test-interval snapshots for the lid-driven cavity problem using 24 POD modes. The first three columns show the reference, lifted OpInf, and original OpInf velocity fields. Color indicates velocity magnitude, while white streamlines and arrows indicate instantaneous flow direction. The last two columns show the corresponding pointwise absolute velocity errors, with insets highlighting the near-lid region.}
    \label{fig:cavity_visualization}
\end{figure}

Finally, Fig.~\ref{fig:cavity_visualization} visualizes the predicted flow fields using 24 POD modes and the corresponding error at five consecutive integer times. While both formulations accurately capture the intricate vortical flow structures, the lifted method achieves lower error, particularly near the top boundary where the control is applied as expected from the guaranteed BC satisfaction of the lifted formulation.

\section{Discussion}
\label{sec:discussion}

The proposed lifted OpInf formulation provides a boundary-consistent treatment of time-dependent nonhomogeneous BCs. By decomposing the solution into a lifting component and a low-rank homogeneous component, the reconstructed state satisfies the prescribed BCs by construction, offering a structural advantage over directly applying OpInf to snapshots with nonzero boundary values. Furthermore, whereas the original OpInf formulation incorporates the control effect through an ad hoc additive linear term, Galerkin projection of the lifted equations naturally yields a richer reduced operator structure that includes bilinear state-input and quadratic input-input interactions for quadratically nonlinear PDEs.

The numerical results demonstrate that the lifted formulation can improve prediction accuracy, projection quality, and solver stability for boundary-actuated systems. At the same time, the lifted and original formulations achieve similar overall accuracy and online computational cost in some settings, supporting the continued use of standard OpInf in many boundary-actuated PDE problems. In particular, when the main objective is accurate prediction of the global state and exact boundary behavior is not critical, applying OpInf with ad hoc linear input terms directly to the original snapshot data may be sufficient. Therefore, the proposed lifted formulation should be viewed not as a necessary replacement for standard OpInf in all boundary-actuated problems, but as a boundary-aware extension that is particularly beneficial when exact satisfaction of the BCs, accuracy near the actuated boundary, or improved stability under boundary forcing is important.

\section{Conclusion}
\label{sec:conclusion}

This paper studied operator inference for nonlinear time-dependent PDEs with inhomogeneous boundary actuation, with the goal of clarifying when boundary lifting is useful in nonintrusive reduced-order modeling. We developed a lifting-based OpInf formulation that decomposes the state into a known boundary-satisfying component and a homogeneous component represented by POD modes. This construction guarantees boundary condition satisfaction by design and leads to a principled reduced model structure involving bilinear state-input and quadratic input-input interactions for PDEs with quadratic nonlinearities.

Numerical experiments on boundary-actuated heat, Burgers, and lid-driven cavity problems showed that the lifted formulation improves boundary compliance and can enhance projection accuracy, prediction accuracy, solver stability, and online efficiency. At the same time, it was observed that standard OpInf can remain competitive in some regimes, especially when the main objective is global state prediction rather than exact boundary consistency. Overall, this study clarifies what lifting contributes to OpInf for boundary-actuated PDEs and when standard OpInf is justified. Future work includes broader comparative studies under more settings and PDEs, parametric boundary controls, noisy or partially observed data, and real-time control or optimization settings.

\section*{Declaration of generative AI and AI-assisted technologies in the manuscript preparation process}

The authors used ChatGPT in a limited capacity to assist with code development and polish writing. The authors reviewed and verified all AI-assisted material and remain
fully responsible for all content of the manuscript.

\section*{Acknowledgments}

Zhuoyuan Wang and Yorie Nakahira were supported in part by the National Science Foundation (NSF) under Grant No.~ECCS-2442948 (CAREER) and in part by the Japan Science and Technology Agency (JST) under Grant No.~JPMJPR2516 (PRESTO). Saviz Mowlavi and Wael Ali were supported solely by MERL.

\section*{Data Availability}

Code and data will be made available upon publication.

\appendix

\section{Additional Numerical Results}

\begin{figure}[t]
    \centering
    
    \begin{subfigure}{0.49\linewidth}
        \centering
        \includegraphics[width=\linewidth]{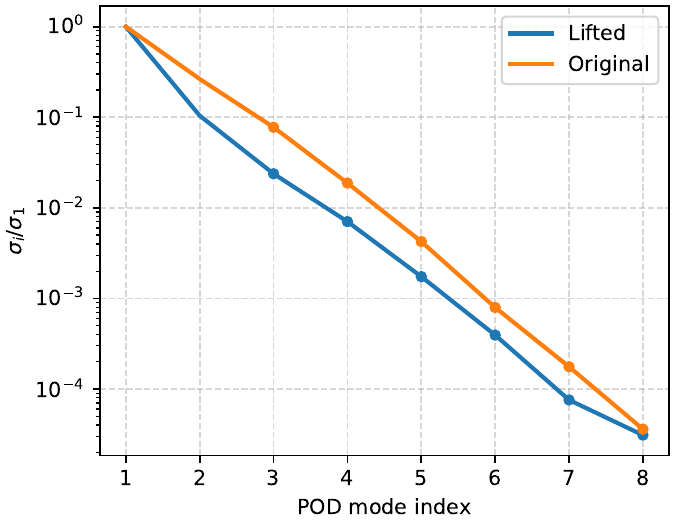}
        \caption{Normalized singular values.}
        \label{fig:heat_pod_singular_values}
    \end{subfigure}
    \hfill
    \begin{subfigure}{0.49\linewidth}
        \centering
        \includegraphics[width=\linewidth]{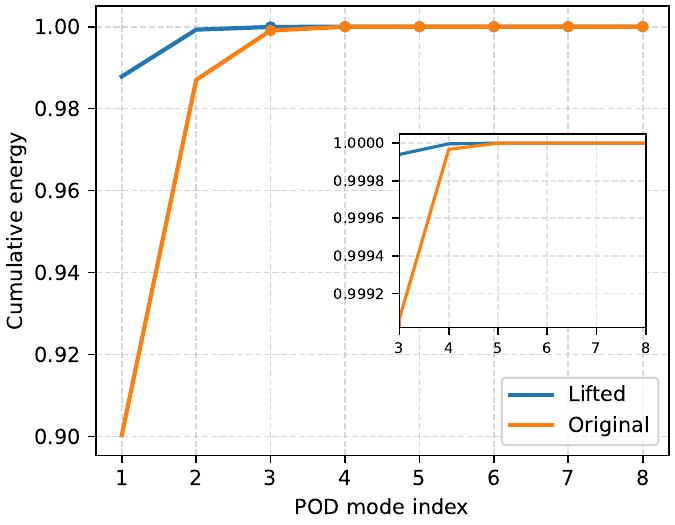}
        \caption{Cumulative captured energy.}
        \label{fig:heat_pod_energy}
    \end{subfigure}
    
    \caption{POD spectrum for the heat equation.}
    \label{fig:heat_spectrum}
\end{figure}

\begin{figure}[t]
    \centering
    
    \begin{subfigure}{0.49\linewidth}
        \centering
        \includegraphics[width=\linewidth]{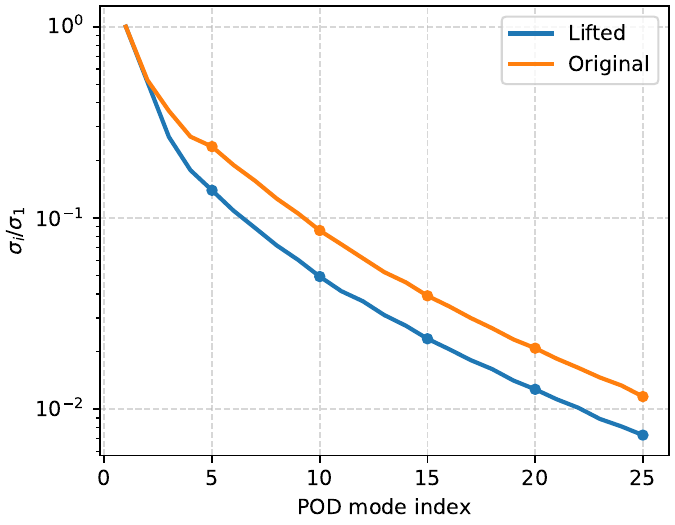}
        \caption{Normalized singular values.}
        \label{fig:burgers_pod_singular_values}
    \end{subfigure}
    \hfill
    \begin{subfigure}{0.49\linewidth}
        \centering
        \includegraphics[width=\linewidth]{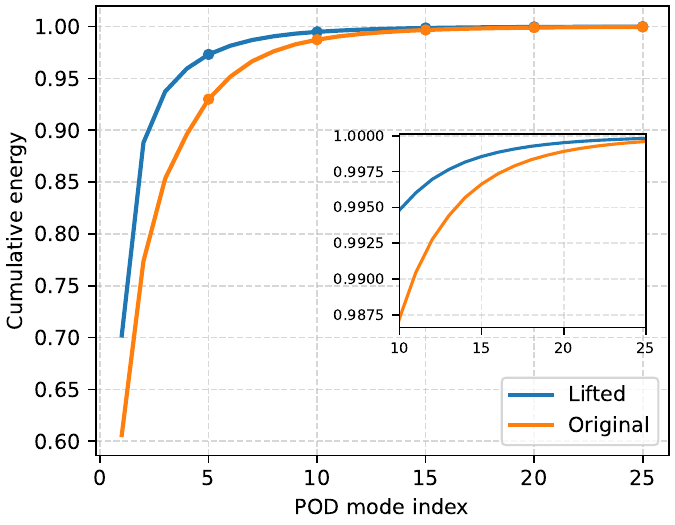}
        \caption{Cumulative captured energy.}
        \label{fig:burgers_pod_energy}
    \end{subfigure}
    
    \caption{POD spectrum for the Burgers' equation.}
    \label{fig:burgers_spectrum}
\end{figure}

\begin{figure}[t]
    \centering
    
    \begin{subfigure}{0.49\linewidth}
        \centering
        \includegraphics[width=\linewidth]{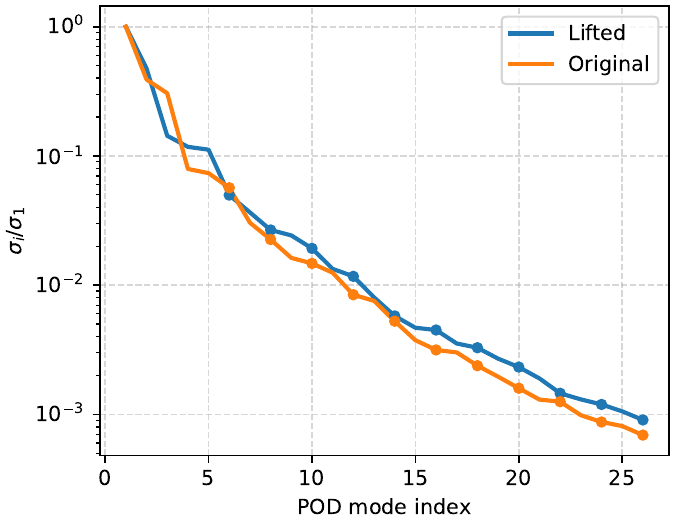}
        \caption{Normalized singular values.}
        \label{fig:cavity_pod_singular_values}
    \end{subfigure}
    \hfill
    \begin{subfigure}{0.49\linewidth}
        \centering
        \includegraphics[width=\linewidth]{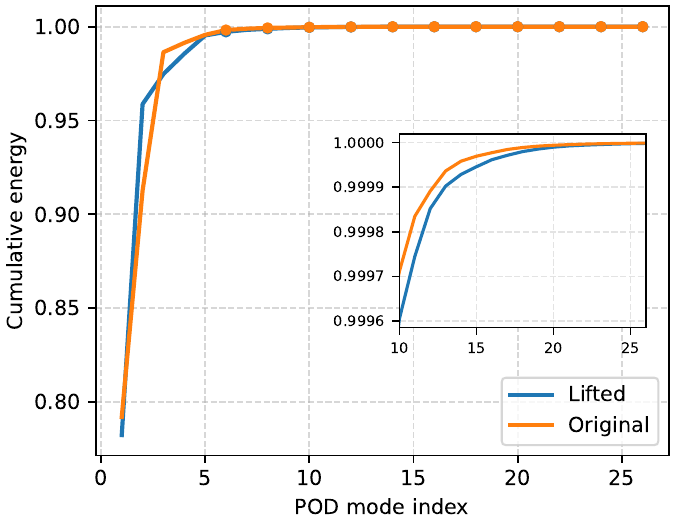}
        \caption{Cumulative captured energy.}
        \label{fig:cavity_pod_energy}
    \end{subfigure}
    
    \caption{POD spectrum for the lid-driven cavity problem.}
    \label{fig:cavity_spectrum}
\end{figure}

Figs.~\ref{fig:heat_spectrum}--\ref{fig:cavity_spectrum} compare the POD spectra of the lifted and original snapshot representations for the heat equation, Burgers' equation, and the lid-driven cavity problem studied in Section~\ref{sec:experiments}. The markers indicate the numbers of POD modes evaluated in the experiments. For the lifted formulation, the POD basis is constructed from the homogeneous component obtained after removing the nonhomogeneous boundary contribution. For the heat and Burgers' examples, the lifted formulation exhibits faster singular value decay and captures a larger fraction of the snapshot energy with fewer POD modes. This suggests that removing the nonhomogeneous boundary contribution before constructing the reduced basis leads to a more compact representation of the remaining homogeneous dynamics, which may contribute to the improved projection and prediction accuracy observed in the main experiments. For the lid-driven cavity problem, the lifted and original spectra are more similar, indicating that the observed accuracy gain of the lifted formulation is not primarily due to better global POD energy capture. Instead, the benefit in this case is more closely related to the boundary-consistent structure of the lifted formulation, which enforces the prescribed boundary conditions by construction and improves the learned reduced dynamics under boundary control.

\begin{figure}[t]
    \centering
    \begin{subfigure}{0.49\linewidth}
        \centering
        \includegraphics[width=\linewidth]{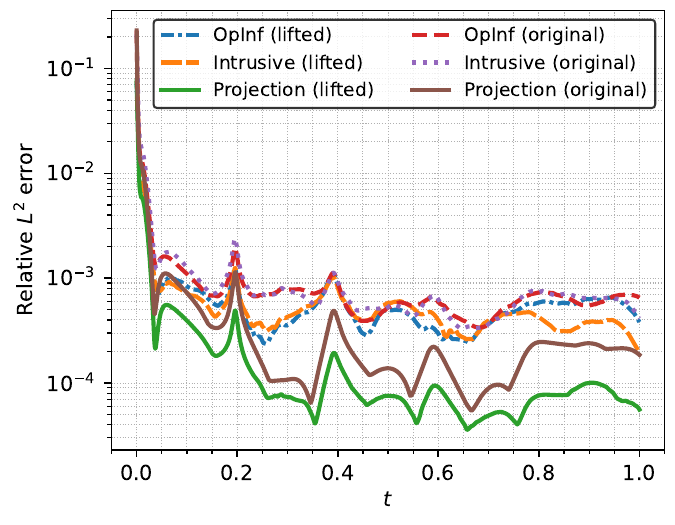}
        \caption{Train}
        \label{fig:heat_train_error_over_time}
    \end{subfigure}
    \hfill
    \begin{subfigure}{0.49\linewidth}
        \centering
        \includegraphics[width=\linewidth]{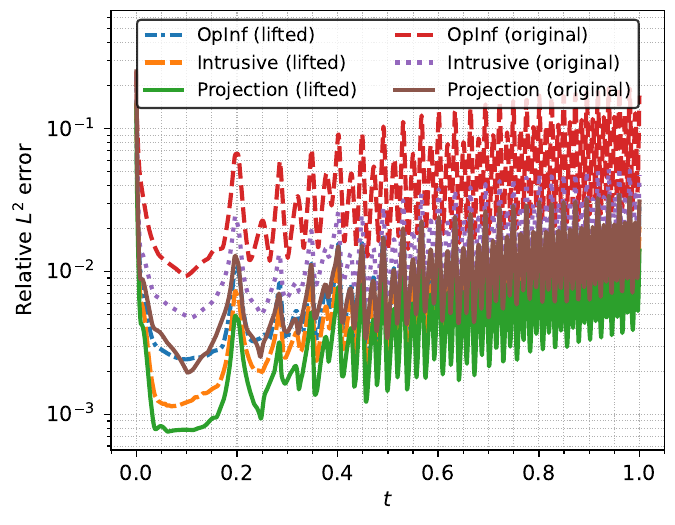}
        \caption{Test}
        \label{fig:heat_test_error_over_time}
    \end{subfigure}
    \caption{Relative $L^2$ prediction error over time for the heat equation, averaged over POD dimensions and successful train-test input combinations.}
    \label{fig:heat_error_over_time}
\end{figure}

\begin{figure}[t]
    \centering
    \begin{subfigure}{0.49\linewidth}
        \centering
        \includegraphics[width=\linewidth]{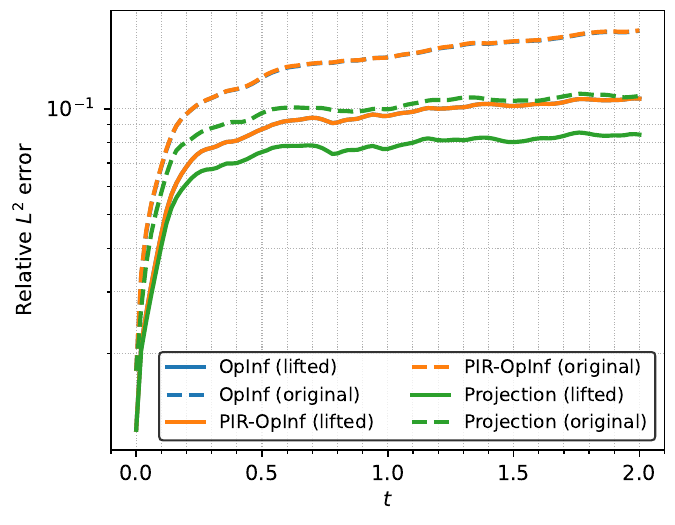}
        \caption{Train}
        \label{fig:burgers_train_error_over_time}
    \end{subfigure}
    \hfill
    \begin{subfigure}{0.49\linewidth}
        \centering
        \includegraphics[width=\linewidth]{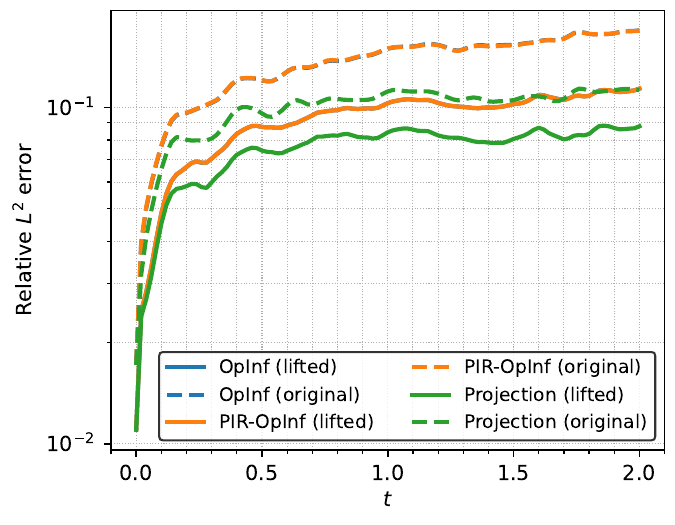}
        \caption{Test}
        \label{fig:burgers_test_error_over_time}
    \end{subfigure}
    \caption{Relative $L^2$ prediction error over time for Burgers' equation, averaged over POD dimensions. For both the lifted and original formulations, the OpInf and PIR-OpInf curves are largely overlapping.}
    \label{fig:burgers_error_over_time}
\end{figure}

\begin{figure}[t]
    \centering
    \includegraphics[width=0.95\linewidth]{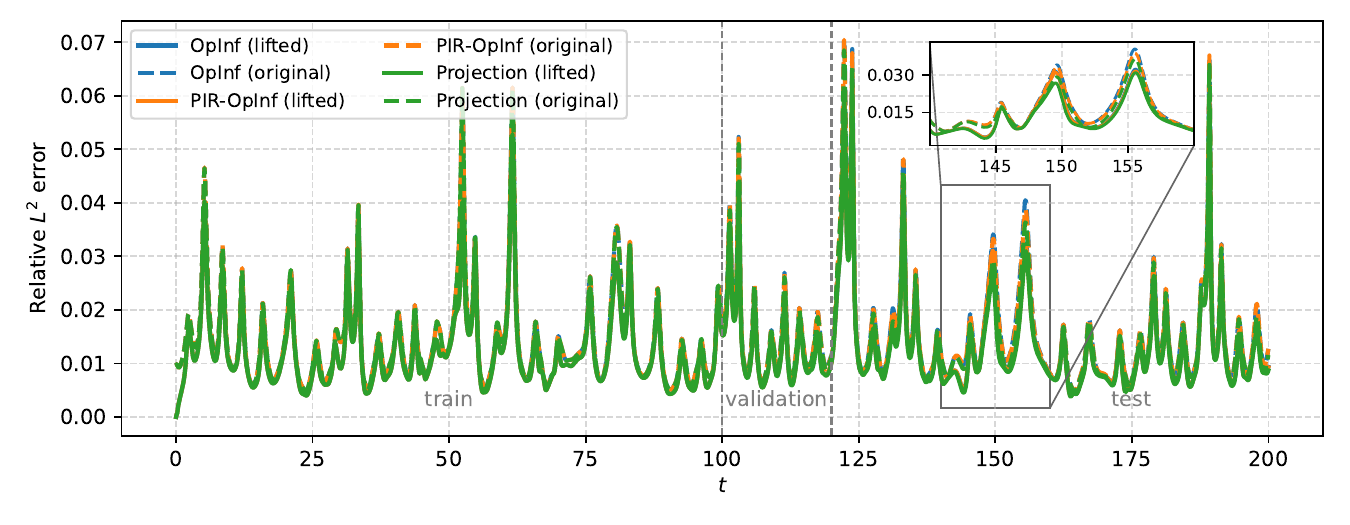}
    \caption{Relative $L^2$ prediction error over time for the lid-driven cavity problem, averaged over POD dimensions. For both the lifted and original formulations, the OpInf, PIR-OpInf, and projection curves are largely overlapping.}
    \label{fig:cavity_error_over_time}
\end{figure}

Figs.~\ref{fig:heat_error_over_time}--\ref{fig:cavity_error_over_time} further show the relative $L^2$ prediction error over time for the heat equation, Burgers' equation, and the lid-driven cavity problem, averaged over POD dimensions. For the heat equation, the averages are computed over successful train-test input combinations only, due to solver failures in the original formulation. Across the three examples, the lifted formulation reduces the relative $L^2$ error over time for all models, including projection, OpInf, PIR-OpInf, and the intrusive ROM, indicating that the benefits of lifting persist throughout the simulation horizon rather than only appearing in time-averaged error values.

\bibliographystyle{elsarticle-num}
\bibliography{citation}

\end{document}